\documentclass[a4paper,10pt]{article}

\usepackage{amssymb,amsmath,amsthm}
\usepackage{braket}
\usepackage{graphicx}
\usepackage{mathrsfs}
\usepackage{subfig}
\usepackage{wrapfig}
\usepackage{amsfonts}
\usepackage[english]{babel}
\usepackage{hyperref}
\usepackage{verbatim} 
\usepackage{standalone}
\usepackage{amsmath}
\usepackage{geometry}
\usepackage{color}
\newcommand{\red}[1]{{\textcolor{red}{#1}}}

\numberwithin{equation}{section}
\usepackage{enumitem}
\usepackage{relsize}
\usepackage{comment}
\usepackage{xcolor}

\input ArtNouvc.fd

\usepackage{aurical}
\usepackage[T1]{fontenc}
\usepackage{pbsi}
\usepackage[T1]{fontenc}
\usepackage{calligra}
\usepackage[T1]{fontenc}

\usepackage{numprint}

\newcommand{\cF}{{\mathcal F}}

\newcommand{\facciatabianca}{\newpage\shipout\null}

\newcommand{\e}{\varepsilon}

\newcommand{\R}{\mathbb R}

\newcommand{\T}{\mathbb T}
\newcommand{\C}{\mathbb C}
\newcommand{\csi}{\xi}

\usepackage{authblk}
\usepackage{float}
\usepackage{amssymb,amsmath,amsfonts,mathrsfs,upgreek,textgreek,bbold,relsize} 
\usepackage{mathabx}
\usepackage{booktabs}
\usepackage{placeins}
\usepackage[nottoc]{tocbibind}
\usepackage{graphicx}
\usepackage{subcaption}

\usepackage{comment}

\renewcommand{\em}{\sl}
\renewcommand\comment[1]{{\iffalse #1 \fi}}

\newcommand\sopra[2]{\genfrac{}{}{0pt}{}{#1}{#2}}

\newtheorem{assumption}{Assumptions}[section]

\newtheorem{remark}{Remark}[section]
\newtheorem{notation}{Notation}[section]
\newtheorem{definition}{Definition}[section]

\newcommand{\io}{{\infty}}

\newcommand\beq[1]{ \begin{equation}\label{#1} }
\newcommand{\eeq}{ \end{equation} }

\newcommand\beqa[1]{ \begin{eqnarray} \label{#1}}
\newcommand{\eeqa}{ \end{eqnarray} }
\newcommand{\beqano}{ \begin{eqnarray*} }
\newcommand{\eeqano}{ \end{eqnarray*} }

\newcommand\dfn[1]{ \begin{definition}\label{#1} }
\newcommand\edfn{ \end{definition} }
\newcommand\ass[1]{ \begin{assumption}\label{#1} }
\newcommand\eass{ \end{assumption} }

\newcommand\notat[1]{ \begin{notation} \label{#1} 
 }
\newcommand\enotat{\end{notation}}

\newcommand\rem{\begin{remark} 
\rm 
}
\newcommand\erem{\end{remark} 
}

\newcommand\equ[1]{{\rm (\ref{#1})}}
\newcommand{\nl}{{\smallskip\noindent}}

\newcommand{\Z}{\mathbb{Z}}
\newcommand{\N}{\mathbb{N}}

\newcommand{\Nf}{\mathtt N}

\newcommand\ttx{{\mathtt x}}
\newcommand\tty{{\mathtt y}}

\def\R{\mathbb R}
\def\T{\mathbb T}

\newcommand\chr\ro

\newcommand{\ro}{{\mathtt r}}

\usepackage{appendix}

\newcommand\ddi\updelta			

\makeatletter
\newcommand{\pushright}[1]{\ifmeasuring@#1\else\omit\hfill$\displaystyle#1$\fi\ignorespaces}
\newcommand{\pushleft}[1]{\ifmeasuring@#1\else\omit$\displaystyle#1$\hfill\fi\ignorespaces}
\makeatother

\begin{document}

\title{{\bf On the analytic properties of the perturbing function in the PCR3Body Problem.}}

\date{14 August 2025}

\author[1]{ Corrado Falcolini} \author[2]{Davide Zaccaria}
\affil[1]{\small Dipartimento di Matematica e Fisica, Universit\`a degli Studi Roma Tre, Roma (Italy)}
\affil[2]{\small Department of Mathematics, University of Toronto, Toronto, Ontario (Canada).}

\maketitle

\noindent
{\bf Abstract.}
We provide a new expansion of the Fourier coefficient of the Perturbing function of the PCR3Body problem in terms of Hansen Coefficients. This gives us a precise asymptotic formula for the coefficient in the region of application of KAM theory (i.e small value of eccentricity and semi-major axis see e.g. \cite{Celletti-Chierchia}). Moreover, in the above region, we study the presence of zeros of the Fourier coefficient for coprime modes $(m,k) \in \Z^2$ and the presence of common zeros between coefficients relative to modes $(m,k)$,$(2m,2k)$ and $(m,k)$,$(2m,2k)$,$(3m,3k)$. Thanks to the previous expansion, this numerical analysis is done up to order $60$ in the power of eccentricity and semimajor axis.  This is a first step for a possible application of \cite{Singular KAM, BBCZ} to PCR3Body Problem that would imply a reduction in terms of measure in the phase space of the so called "non--torus" set from $O(1-\sqrt{\e})$ (implied by  standard KAM theory) to $O(1-\e |\log\e|^c )$ for some $c>0$.

 Planar circular restricted three-body problem. Zeros of Fourier coefficients.


\tableofcontents

\section{Introduction}\label{introduction}

The Kolmogorov–Arnold–Moser (KAM) theory addresses the rigorous construction of quasi-periodic trajectories in nearly--integrable Hamiltonian systems. To be more precise, a classical result in this theory can be stated as follows: {\it consider a real-analytic Hamiltonian system $H_\e(x,y)=h(y) + \e f(x,y)$ defined on the phase space $\T^n \times B \subset \R^{2n}$. Let $y_0 \in B $ such that $\omega(y_0)=h'(y_0)$ is a Diophantine vector \footnote{i.e. $|\omega(y_0) \cdot k|>\gamma |k|^{-\tau} \ \forall \ k \in \Z^n\setminus \{0\}$ with $\gamma,\tau>0$} and that $\det h''(y_0) \ne0$. Then there exists $\e_0$ small enough such that for any $0<\e<\e_0$ the unpertubed torus $\T^n \times \{y_0\}$ can be analytically conjugated to an invariant torus on which the flow of $H_\e$ is analytically conjugated to a linear flow with frequency $\omega(y_0)$.}
This formulation implies quite easily  (see e.g. \cite{sqrte,Psqrt}) that the preserved KAM invariant tori form a set of positive Liouville (Lebesque) measure, whose complement has a measure proportional to $\sqrt{\e}$. In this setting, in which we take into account only the so-called "primary" tori (those preserved by the perturbation), this measure estimate is {\sl optimal} in the sense that there appear $O(\sqrt{\e})$ region near simple resonances free of primary invariant tori (this can be seen for example with the trivial 1 degree of freedom system $p^2/2 + \e \cos q$).

\nl
Given the above situation, in order to improve the estimates on the total measure of the region filled with invariant tori, it is clear that one has to carefully study the existence of the so--called "secondary" tori, namely quasi-periodic motions that were not present in the unperturbed case and were formed by the presence of the perturbation (think also in this case to the phase space enclosed by separatrix in the 1-dimensional $p^2/2 + \e \cos q$). The rigorous analysis of such motions is remarkably complicated and has been recently done by Biasco-Chiercia in a sequel of works \cite{Complex Arnold,NL,GP,Singular KAM} for natural systems (i.e. $y^2 + \e f(x)$) and then has been extended to convex systems $h(y) + \e f(x)$ in \cite{BBCZ}. For these kind of systems, under some generic condition on the perturbation, they improved the previous estimates, bringing to a set of initial data that lie in an invariant torus with volume of order $\e |\log \e|^c$ for some $c>0$ in agreement (up to logarithmic correction) to a Arnol'd--Kozlov--Neishtadt conjecture in \cite{AKNDUE}.

\nl
As it is well known, the development of KAM theory was inspired by classical challenges in Celestial Mechanics, as the {\it n}-body problem, and so a crucial and long-standing problem in this area is the application of such theory to physical systems for "experimental" values of the involved parameters. In this direction, the first difficulty that arises is non--degeneracy of the integrable part. \footnote{For example, typical n-body problems are strongly degenerate (see e.g. \cite{Herman}).} Fortunately, there is a model in celestial mechanics that is KAM non-degenerate and has therefore attracted a lot of attention and became a "special" model for the application of KAM theory. This model is the {\it planar, circular, restricted three-body problem} (PCR3Body Problem or simply PCR3BP).\footnote{PCR3BP is actually isoenergetically non--degenerate, i.e. $\det \begin{pmatrix}
    h''(y_0) & \omega(y_0) \\ \omega(y_0) & 0
\end{pmatrix} \ne 0$ but KAM theorem can be easily adapted to those systems (see \cite{Celletti-Chierchia}).}

\nl
The PCR3BP describes the bounded planar motion of a {\it "zero mass"} body subject to the gravitational field generated by two primary bodies, in the same plane, revolving  on circular Keplerian orbits (which are assumed to be not influenced by the small body)

\nl
When the mass ratio of the two primary bodies is small the PCR3BP is described by a nearly-integrable Hamiltonian system with two degrees of freedom that physically represent the eccentricity of the "zero mass" body  $e \in [0,1)$ and its semimajor axis $a \ge 0$.
In \cite{Celletti-Chierchia} classical KAM theory is applied to such model and, since the phase space is four dimensional and the energy level has three dimensions, then the invariant tori have two dimension and so there is no room to escape. Thus we say that this system is {\it totally stable} in the sense that in a
neighborhood of any phase point of negative energy, if
the mass ratio of the primary bodies is small enough,
the asteroid stays forever on a nearly Keplerian ellipse.

\nl
Summing up, PCR3BP is one of the most relevant model for applications in Celestial Mechanics and, moreover, recent development of KAM theory have led to a significant improvement in stability estimates for a very wide class of systems.
So, two questions arise naturally:

\vskip0.3cm

\centering{\textit{Which is the total measure of quasi--periodic motions for PCR3BP?}}

\vskip0.3cm

\centering{\textit{Is it possible to apply the recent theory of \cite{Singular KAM, BBCZ} to PCR3BP improving the estimate given by standard KAM theory?}}

\flushleft
Toward this direction, the first necessary task is to check if the Biasco-Chierchia condition on the perturbation could fit with our model.
The condition (that is proved to be generic and prevalent in the space of real-analytic function see \cite{GP}) reads as follow\footnote{In \cite{Singular KAM,BBCZ} the setting is more general and has $n$ dimension, but since we are dealing with a 2-degree-of-freedom system we will only consider that case.}: let $f: \T^2 \to \R$ be real analytic with radius $s>0$, and denote with $f_{m,k}$ its Fourier coefficient of mode $(m,k) \in \Z^2$, then 
\begin{equation} \label{condition BC}
\varliminf_{\sopra{|m|+|k|\to+\io}{(m,k)\in\mathcal{G}^2}} |f_{m,k}| \, \, e^{(|m|+|k|)s}  \, \, (|m|+|k|)^2>0\,,
\end{equation}
where 
$$
\mathcal{G}^2:=\{(m,k)\in \Z^2\setminus \{0\}: \mbox{the first non–null component is strictly positive and} \ {\rm gcd} (m,k)=1\}\,.
$$ 
This easily implies that there exists $0\le \delta<1$ and a suitable (but $\e$-indepent) $\Nf = \Nf(\delta,s)$, such that
\begin{equation} \label{condition 2}
|f_{m,k}|\geq \delta (|m|+|k|)^{-2}\, e^{-(|m|+|k|) s}\,,\qquad \forall \ (m,k)\in\mathcal{G}^2\,,\ |m|+|k|\ge\Nf.
\end{equation}

\nl
This condition is extremely useful since near simple resonance, after high order averaging (see \cite{NL}), the system $h(y) + \e f(x)$ behaves like
$$
h(y) + \e g(y) + \e \pi_{\Z (m,k)} f (m x_1 + kx_2) + \mbox{exp. small reminder}
$$
for some real-analytic $g$ and where $\pi_{\Z (m,k)} f (\theta) = \sum_{j \in \Z} f_{jm,jk} (\theta) e^{ij \theta}$.

If \ref{condition 2} holds, it can be shown (see e.g. \cite{GP}) that for suitable $\phi_k \in \R$
$$
\e \pi_{\Z (m,k)} f (\theta) \sim  \e |f_{m,k}| \cos (\theta + \phi_k) \ \ \mbox{for} \ \ |m|+|k|>\Nf
$$
and so, near simple resonance, instead of dealing with an infinite number of different arbitrary systems, one has to study an infinite number of approximately 1-dimensional simple pendula (i.e. $\sim h(\tty_1)+ \e |f_{m,k}| \cos(\ttx_1+\phi_k)$ for $\ttx_1=mx_1 + kx_2$) and just a finite number of arbitrary systems (the ones for $|m|+|k|\le \Nf$). In this way, it is possible to have a uniform and quantitative control on analytic properties of action angles variable in order to reach a high precision on the measure estimates on secondary tori.

\nl
Coming back to condition \ref{condition 2}, it clearly implies that
$$
|f_{m,k}| > 0 \ \  \forall \ ({m,k}) \in \mathcal{G}^2, |m|+|k|\ge\Nf
$$
and so the first step for applications is to check the absence of zeros of fourier coefficients for high fourier modes.

\nl
One can notice a crucial difficulty that arises immediately in the application to PCR3BP: in our model the perturbation depends on both angles $x$ and actions $y$ \footnote{see \ref{hamiltoniana} in section \ref{section2}.} and so its Fourier coefficients $f_{m,k}(y)$ is now an holomorphic function that, in general, does have zeros in its domain. In these cases, as suggested by L. Niedermann and S. Barbieri in \cite{TesiSanti}, the condition is replaced by the fact that, roughly speaking, it exists a natural number $\ell \ge 2$ such that
$$
| \sum_{j=1}^\ell f_{(jm,jk)}(y)|>0 \ \  \forall \ (m,k) \in \mathcal{G}^2, \ \ |m|+|k|>\Nf
$$
since, generically and also in a measure-prevalent sense, $\ell$ holomorphic function have essentially not a common zero in their domain (see Sard-Morse theory developed by Yomdin-Comte in \cite{YC} for a quantitative and rigorous statement).

\nl
In this work we do a numerical study of the zeros of the Fourier coefficients $f_{m,k}$ of the perturbing function of PCR3Body problem. In particular we check, up to truncation, if the coefficients relative to coprime modes satisfy the conditions written above for a possible application of KAM theory for secondary invariant tori.

\nl
In order to do this analysis, following \cite{LB10}, in section \ref{section2} we provide an expansion of the $f_{m,k}$ in terms of Hansen coefficients that allow us to reach in a short time a high order of truncation. This is necessary since, as we can see in figure \ref{Truncationdouble} and \ref{truncationtriple}, the curves that describe the zeros change significantly as truncation changes and usually stabilize around order 50/60 in power of $e,a$. 

\nl
In addition to this numerical power, the expansion in terms of Hansen coefficients gives us a sharp asymptotic formula for small values of $(e,a)$ (the KAM regime of application) i.e. for $(m,k) \ne 0$ we can write
$$
f_{m,k} (a,e) = 2 \, t_{m,k} e^{|m-k|} a^{m^*} [1+O(e^2,a)], \qquad m^* = \begin{cases}
    m, \ \ \mbox{if} \ m \ge 2 \\
    m+2  \ \ \mbox{if} \ m \in \{0,1\}.
\end{cases}
$$
and we provide the precise formula for the coefficients $t_{m,k}$ in \ref{tmkA} for $m\ge 2$ and in \ref{tmkB} for $m \in \{0,1\}$.

\nl
In section 3, using the so--called \textit{Wnuk's method} (see \cite{cinesi}), we provide tables with truncated Hansen coefficients $X^{n,m}_k$ for $k=0,1,4,8$. Using those tables and the previous expansion one can obtain a truncated formula for Fourier coefficient of PCR3BP up to high order in power of $a,e$. This is shown in appendix \ref{appendicite}.

\nl
With the numerical power developed for those coefficients, we can eventually study the curves described by $f_{m,k}(a,e)=0, (a,e) \in (0,1)\times (0,1)$ and $(m,k) \in \mathcal{G}^2$. Plotting those curves we can see in figure \ref{ZERISINGOLI} that, increasing the order of truncation in powers of $a$ and $e$, the number of curves (i.e. of zeros) increases and their distance from the origin decline towards zero. This implies that is not possible to fix a region of small values $(a,e)$ such that $|f_{m,k}(a,e)|>0$ uniformly at every order. Since this first property fails, we pass to see the points for which simultaneously $f_{m,k} (a,e) =0$ and $f_{2m,2k}(a,e)=0$. The picture shown in figure \ref{ZERIDOPPI} gives the same result as the previous case: most of the double zeros are far from the origin but few of them come closer and closer to the origin while the order of truncation increases.
Finally, a different result holds in the last case in which we check the presence of points $(a,e)$ such that $|f_{m,k} (a,e)|+|f_{2m,2k} (a,e)|+|f_{3m,3k} (a,e)|=0$. Indeed, we do not find any point satisfying this condition up to order 60, but we show that there is a mode $(m,k)=(5,-2)$ and a point $(a,e)=(0.188,0.900)$ at which the distance from those three curves is of order $10^{-5}$ (see figure \ref{zerizoom52}). Moreover, in this "critical case", high fourier modes appear so that the curves of zeros are not stable at order 60. So we cannot reject the hypothesis of a common zero among the three coefficients at a higher order of truncation. A more refined analysis would be needed.


\section{Expansion and asymptotic of the Fourier coefficients} \label{section2}
Let $P_0,P_1$ and $P_2$ be three bodies of masses, respectively, $m_0,m_1$ and $m_2$. We assume that $m_2$ is much smaller than $m_0$ and $m_1$ (restricted problem) and that the motion of $P_1$ around $P_0$ is circular. We also assume that the three bodies always move on the same plane.
Let $\mu >0$ be defined as $\mu := m_0^{-2/3}$ so that the two-body Hamiltonian in Delaunay action--angle variables $(\Lambda,\Gamma,\lambda,\gamma)$ \footnote{By Kepler's laws one finds that $\Lambda = \mu \sqrt{(m_0+m_1) \, a}$, and $\Gamma = \Lambda \sqrt{1-e^2}$ where $a$ is the semimajor axis and $e$ is the eccentricity of the Keplerian orbit of $P_1$ around $P_0$. } (see \cite{Brouwer}, \cite{Celletti-Chierchia}, \cite{CM}, \cite{VS}, etc.) becomes
$$H_D (\Lambda,\Gamma,\lambda,\gamma) = - \frac{1}{2 \Lambda^2}$$
while we introduce the perturbing parameter as $\e= \frac{m_1}{m_0^{2/3}}$.

\nl 
Set the units of measure so that the distance between $P_0$ and  and $P_1$ is one and so that $m_0+m_1=1$.
Taking into account the interaction of $P_2$ on $P_1$, the Hamiltonian function governing the PCR3BP becomes
\begin{equation} 
    H(\Lambda,\Gamma,\lambda,\gamma,t) = - \frac{1}{2 \Lambda^2} + \e (r_1 \cos{(\varphi - t) - \frac{1}{\sqrt{1+r_1^2 - 2 r_1 \cos{(\varphi - t)}}}}\,,
\end{equation}
where $r_1$ is the distance between $P_0$ and $P_1$ and $\varphi$ is the longitude of $P_1$ w.r.t. $P_0$.
Then, $\varphi= f + \gamma$ where $f$ is called {\it true anomaly}. Finally, by the generating function of Delaunay--variables, we know that the {\it mean anomaly} $\lambda$ is an implicit function of the {\it eccentric anomaly} $u$ (i.e.  the angle between the origin of the plane and the planet during its motion) and eccentricity $e = e(\Lambda,\Gamma)=\sqrt{1 - ( \frac{\Gamma}{\Lambda})^2}$ by the relation
\begin{equation}
    \lambda = u - e\, \sin{u}.
\end{equation}
Using these information, with the canonical change of variables
\begin{equation} \label{delaunayvar}
    L=\Lambda, \quad \ell=\lambda; \qquad G= \Gamma, \quad g=\gamma-t; \qquad\qquad (L,G,\ell,g) \in \{ 0 < G <L\} \times \T^2;,
\end{equation}
the Hamiltonian becomes
\begin{equation} \label{hamiltoniana}
H_{rpc} ( L, G,\ell, g)= - \frac{1}{2\,L^2} - G + \e F(L, G,\ell, g)
\end{equation}
where $F(L, G,\ell, g)$ is the perturbing function
\begin{equation} \label{perturbation}
    R(r_1,\varphi,t) = r_1 \cos(\varphi - t) - \frac{1}{\sqrt{1 + r_1^2 - 2 r_1 \cos(\varphi - t)}}  
\end{equation}
but with $r_1$ and $\varphi - t$ expressed in terms of the Delaunay variables $(L, G,\ell, g)$.
Now we need to express the perturbing function in \equ{perturbation} in its Fourier series studying its Fourier coefficients.
For this kind of function that arises from gravitational force in $1782$ Legendre finds an useful expansion on some polynomials (then extended by Laplace) such that they can be expressed by the generating formula
\begin{equation} \label{generatrice}
    \frac{1}{\sqrt{1+z^2-2z \csi}} = \sum_{n=0}^\infty z^n P_n (\csi); \quad (\mbox{in our case} \ \csi = \cos(\varphi - t), z=r_1)
\end{equation}
These Legendre polynomials that we have called $P_n$ can be expressed as
\begin{equation}
    P_n (\csi)= \sum_{k=0}^{[n/2]} p_{n,k} \csi^{n-2k}
\end{equation}
where the notation $[n/2]$ represent the largest integer less than or equal to $n/2$ (namely, the floor of $n/2$) and
\begin{equation}
    p_{n,k} = \frac{(-1)^k}{2^n} \frac{(2n - 2k)!}{k! (n-k)! (n-2k)!}.
\end{equation}
In this way one can obtain the following expasion 
\begin{equation}
    \begin{split}
    R(r_1,\varphi,t) &= r_1 \cos(\varphi - t) - P_0(\cos(\varphi - t)) - r_1 P_1(\cos(\varphi - t)) - \sum_{j=2}^\infty P_j (\cos(\varphi-t)) \,r_1^j
    \\
    &= - (1 + \sum_{j=2}^\infty P_j (\cos(\varphi-t))\, r_1^j ).
    \end{split}
\end{equation}
In order to obtain some better expansion, one can do a classical computation (see \cite{LB10},\cite{WW27}) that is for all $r_1 \in [0,1)$
\begin{equation}
    \begin{split}
        (1 + r_1^2 - 2 r_1 \cos(\varphi - t))^{-1/2} &= \big(1 -  r_1 \, \mbox{exp}\,(-i(\varphi - t)) \big)^{-1/2} \big(1 -  r_1 \mbox{exp}\,(-i(\varphi - t))\big)^{-1/2} 
        \\
        &= \sum_{p=0}^\infty \sum_{q=0}^\infty \frac{(2p)!(2q)!}{2^{2p+2q} (p!)^2 (q!)^2} \mbox{exp}\,\big(i(p-q)(\varphi - t)\big) r_1^{p+q}
    \end{split}
\end{equation} 
that, with $n=p+q$, and after changing $q$ to $n-q$, becomes
\begin{equation}
    \qquad \qquad \qquad \qquad \qquad \qquad  = \sum_{n=0}^\infty \bigg( \sum_{q=0}^n \frac{(2q)! (2n-2q)!}{2^{2n} (q!)^2 ((n-q)!)^2} \mbox{exp}\, \big(i(2q-n)(\varphi - t)\big)  \bigg)r_1^{n}.
\end{equation}
Thus, comparing two expression above with \equ{generatrice}, for $0\le q \le n$ we have
\begin{equation}
    P_n (\cos(\varphi-t)) = \sum_{q=0}^n \widetilde{\cF}_{q,n}\, \mbox{exp}\, \big(i(2q-n)(\varphi - t)\big)  := \sum_{q=0}^n \frac{(2q)! (2n-2q)!}{2^{2n} (q!)^2 ((n-q)!)^2} \mbox{exp}\, \big(i(2q-n)(\varphi - t)\big) 
\end{equation}
such that the perturbing function becomes
\begin{equation}
    \begin{split}
    R &= -1- \sum_{n=2}^\infty \bigg(\frac{r_1}{a}\bigg)^n \bigg(\sum_{q=0}^n \widetilde{\cF}_{q,n} \, \mbox{exp} \,\big({i(2q-n)(\varphi - t)\big)}\bigg) a^n
    \\
    &= -1- \sum_{n=2}^\infty \sum_{q=0}^n {\cF}_{q,n} \,  \mbox{exp} \,\big({i(2q-n)(\varphi - t)\big)} \, a^n; \qquad  {\cF}_{q,n}:=\widetilde{\cF}_{q,n}\bigg(\frac{r_1}{a}\bigg)^n,
    \end{split}
\end{equation}
using Delaunay variables in \equ{delaunayvar} we  write $\varphi-t = f+ g$ and the above expression becomes
\begin{equation}
    F = -1- \sum_{n=2}^\infty \sum_{q=0}^n {\cF}_{q,n} \,  \mbox{exp} \,\big({i(2q-n)f\big)} \,  \mbox{exp} \,\big({i(2q-n)g\big)} \, a^n.
\end{equation}
where $F$ is the perturtubing function in Delaunay variables, i.e. is the function $R$ composed with the change of coordinates in \equ{delaunayvar}.
Now the crucial point is that the coefficients ${\cF}_{q,n}$ can be expressed in terms of Delaunay variables $(\ell,g,a(L),e(L,G))$ via the {\it Hansen coefficients} $X^{n,m}_k (e)$ (section \ref{HANSEN}) defined for $n,m \in \Z$ such that
\begin{equation} \label{Giacaglia}
    \bigg(\frac{r_1}{a}\bigg)^n \,  \mbox{exp} \,\big({imf\big)} = \sum_{k=-\infty}^{+\infty} X^{n,m}_k (e) \ \mbox{exp} \,\big({ik \ell\big)}.
\end{equation} 
\noindent Thus we have
\begin{equation}
    \begin{split}
    F &= -1- \sum_{n=2}^\infty \sum_{q=0}^n  \sum_{k=-\infty}^{+\infty} \widetilde{\cF}_{q,n} X^{n,2q-n}_k (e) \,  \mbox{exp} \,\big({ i[(2q-n)g + k \ell]\big)} \, a^n.
    \\
    &= -1- \sum_{n=2}^\infty \sum_{m=-n}^{n}{\vphantom{\sum}}'  \sum_{k=-\infty}^{+\infty}  \widetilde{\cF}_{\frac{m+n}{2},n} X^{n,m}_k (e)  \,  \mbox{exp} \,\big({ i(mg + k \ell)\big)} \, a^n
    \\
    &=  -1- \sum_{n=2}^\infty \sum_{m=-n}^{n}{\vphantom{\sum}}'  \sum_{k=-\infty}^{+\infty} \frac{(m+n)! (n-m)!}{2^{2n} ((\frac{m+n}{2})!)^2 ((\frac{n-m}{2})!)^2} X^{n,m}_k (e)  \,  \mbox{exp} \,\big({ i(mg + k \ell)\big)} \, a^n
    \end{split}
\end{equation}
where $\sum{\vphantom{\sum}}'$ indicates that the sum is over every term between $-n$ and $n$ separated by $2$.

\nl
Thanks to the parity of the perturbing function, we just know that
\begin{equation} 
\begin{aligned}
\label{Espansione r,k}
    F (a,e,\ell,g) &= -1- \sum_{n=2}^\infty \sum_{m=-n}^{n}{\vphantom{\sum}}'  \sum_{k=-\infty}^{+\infty}  \frac{ (m+n)! (n-m)!}{2^{2n} ((\frac{m+n}{2})!)^2 ((\frac{n-m}{2})!)^2} X^{n,m}_k (e)  a^n \cos (mg + k \ell) \\
    &:= -1- \sum_{n=2}^\infty \sum_{m=-n}^{n}{\vphantom{\sum}}'  \sum_{k=-\infty}^{+\infty}  C_{n,m} X^{n,m}_k (e)  a^n \cos (mg + k \ell)
\end{aligned}
\end{equation}
so that, for fixed $(m,k)$ in $\N \times \Z$ one can obtain the expression for the { \it Fourier coefficient of the perturbing function} as
\begin{equation} \label{HansenExpansion}
    f_{m,k}(a,e) = 
    \begin{cases} \mbox{if} \ \ (m,k)=(0,0)  \quad 
    -1 -  \displaystyle\sum_{n=2}^\infty{\vphantom{\sum}}' \ C_{n,0} X^{n,0}_0 (e)\,  a^n, \\
   \mbox{if} \ \ (m,k)\neq(0,0) \quad -2\displaystyle\sum_{n=m^*}^\infty\vphantom{\sum}'  \ C_{n,m} X^{n,m}_k (e)\,  a^n 
\end{cases}
\end{equation}
where $m^*=m+2$ if $m=0$ or $m=1$ and $m^*=m$ otherwise. The presence of $m^*$ is due to the fact that in \ref{Espansione r,k} $n$ starts from $2$, while the multiplicative factor $2$ in front of the second sum is due to the parity of cosine.

\nl
Thanks to \cite{Bal05} we can use the following expansion for Hansen coefficients for $k=m+s$ and in the case of $s=k-m>0$:
\begin{equation}\label{balmino}
\begin{split}
    X^{n,m}_{m+s}(e) = &(-1)^s \bigg( \frac{e}{2} \bigg)^s \sum_{t=0}^{\infty} \bigg\{ \sum_{j=0}^t \sum_{p=0}^j \binom{n+m+1}{j-p} \frac{(m+s)^p}{p!} \sum_{q=0}^{s+j} \binom{n-m+1}{s+j-q} \frac{(m+s)^q}{q!} (-1)^q
    \\
    &\bigg[ 2 \binom{2t -n +s -p-q-2}{t-j} - \binom{2t -n +s -p-q-1}{t-j}   \bigg]  \bigg\} \bigg( \frac{e}{2} \bigg)^{2t}
\end{split}
\end{equation}
where binomial coefficient $\binom{-\mu}{p}$ where $\mu>0$ must be computed as being equal to $(-1)^p \binom{\mu + p -1}{p}$, $p$ being always positive.

\nl 
If $s=k-m<0$, one can use the fact that $X_{k}^{n,m} = X_{-k}^{n,-m}$ to calculate $ X^{n,m}_{m+s}(e) = X^{n,-m}_{-m-s}(e)$ using the formula in \ref{balmino}.
In this way it is easy to see that, for $s=k-m \ge 0$, the leading term is for $t=0$, i.e. it is $e^{k-m}$, and the same holds for $-s=m-k \ge 0$ with the right change of sign, namely that the leading term is $e^{m-k}$. So it is clear that $X^{n,m}_{k}(e) = o(e^{|k-m|})$, i.e. using \ref{HansenExpansion}, that
\begin{equation}
\label{ASINTOTICA}
    f_{m,k}(a,e)=
    \displaystyle\begin{cases}
        \mbox{if} \ (m,k)=(0,0) \quad -1+ t_{0,0}\,  \,a^{2} \, \big[ 1 + \mathcal{O}(e^2;a) \big] \\[5pt]
        \mbox{if} \ (m,k)\neq(0,0) \quad
    2\,t_{m,k}\, e^{|m-k|} \,a^{m^*} \, \big[ 1 + \mathcal{O}(e^2;a) \big] 
    \end{cases}
\end{equation}
where $m^*=m+2$ for $m=0,1$ and $m^*=m$ otherwise, and so the form of the coefficient $t_{m,k}$ control the asymptotic behaviour of the Fourier coefficient for small values of $(a,e)$. Since in expansion \ref{HansenExpansion} we can see that the first term of the expansion is for $n=m^*$, we have to divide our analysis in the case in which $m^*=m$ and $m^*=m+2$. We start with the easiest case.

\nl
{\bf Case A:} ${\mathbf{m \ge 2 \ \   \mbox{i.e.} \ \  m^*=m.}}$

\nl
From expansions \ref{HansenExpansion} in this case we know that
\begin{equation}
    t_{m,k}^{(A)} = -\frac{(2m)!}{2^{2m} (m!)^2} [X^{m,m}_k]
\end{equation}
where $[X^{m,m}_k]$ indicates the coefficient that multiplies the term $e^{|k-m|}$ in the $e$-power expansion of the Hansen coefficient.

\nl
As we can see from \ref{balmino}, in the series expansion of the Hansen coefficient, the eccentricity has a power expressed by $s+2t$, so in order to evaluate the coefficient of the leading term $e^{|k-m|}$, we have to look at the term with
$$
s+2t=|k-m|
$$
so it is clear that we have to distinguish between two cases:

\nl
Case A.1: ${\mathbf{m-k<0}}$:
for our intent, w.r.t. \ref{balmino} we set $n=m$ and $s=k-m$, and if we want the term $t_{m,k}$ we need $s+2t=|k-m|$ that means, for $m<k$, $k-m+2t=k-m \Rightarrow t=0$.
\begin{equation}
    [X^{m,m}_{k}] = \frac{(-1)^{k-m}}{2^{k-m}}  \binom{2m+1}{0}  \sum_{q=0}^{k-m} \binom{1}{k-m-q} \frac{k^q}{q!} (-1)^q \bigg[ 2 \binom{k -2m-q-2}{0} - \binom{k -2m-q-1}{0}   \bigg].
\end{equation}
so that the first binomial coefficient that involves $q$ is different from zero only if $q=k-m$ or $q=k-m-1$ and such that
\begin{equation}
\begin{split}
    [X^{m,m}_{k}] &= \frac{(-1)^{k-m}}{2^{k-m}}\big[\frac{k^{k-m}}{(k-m)!} (-1)^{k-m} - \frac{k^{k-m}}{(k-m)!} \frac{k-m}{k} (-1)^{k-m} \big] =\frac{k^{k-m}}{2^{k-m} \,(k-m)!} [1 -\frac{k-m}{k}]
    \\
    &=  \frac{1}{2^{k-m}} \frac{m \, k^{k-m}}{k\, (k-m)!}.
\end{split}
\end{equation}
Finally, denoting with $t_{m,k}^{(A,-)}$ the value of $t_{m,k}^{(A)}$ in the case $k>m$, we have obtained
\begin{equation}
 t_{m,k}^{(A,-)} =- \frac{(2m)!}{2^{m+k} (m!)^2} \frac{m \, k^{k-m}}{k\, (k-m)!}.
\end{equation}

\nl
Case A.2: ${\mathbf{m-k>0}}$:
Now, setting always $s=k-m$, in order to control only the dominant term, we have to ask that $s+2t=|k-m|=m-k$ namely that $t=m-k$. So from \ref{balmino} we have
\begin{equation}
\begin{split}
    [X^{m,m}_{k}] = &(-1)^{k-m} \bigg( \frac{1}{2} \bigg)^{m-k} \bigg\{ \sum_{j=0}^{m-k} \sum_{p=0}^j \binom{2m+1}{j-p} \frac{k^p}{p!} \sum_{q=0}^{k-m+j} \binom{1}{k-m+j-q} \frac{k^q}{q!} (-1)^q
    \\
    &\bigg[ 2 \binom{ -k-p-q-2}{m-k-j} - \binom{-k-p-q-1}{m-k-j}   \bigg]  \bigg\}.
\end{split}
\end{equation}
We start analyzing the sum in $q$; in order to make the binomial coefficient different from zero, one finds the condition $k-m+j \ge 0$, i.e. $j \ge m-k$, but since the index $j$ goes from $0$ to $m-k$, the only possible contribution is from $j=m-k$, and so $q=0$. In this way the expression becomes
\begin{equation}
\begin{split}
    [X^{m,m}_{k}] &= (-1)^{k-m} \bigg( \frac{1}{2} \bigg)^{m-k} \bigg\{  \sum_{p=0}^{m-k} \binom{2m+1}{m-k-p} \frac{k^p}{p!} \bigg[ 2 \binom{ -k-p-2}{0} - \binom{-k-p-1}{0}   \bigg]  \bigg\}
    \\
    &=\frac{(-1)^{k-m}}{2^{m-k}} \bigg\{  \sum_{p=0}^{m-k} \binom{2m+1}{m-k-p} \frac{k^p}{p!}  \bigg\}.
\end{split}
\end{equation}
Notice that this sum has contributes different from zero until $2m+1 \ge m-k-p$, i.e. $p \ge -(m+k+1)$. So if $m+k+1<0$ the contributions different from zero start from $p=-(m+k+1)$.

\nl
So finally, denoting with $t_{m,k}^{(A,+)}$ the value of $t_{m,k}^{(A)}$ in the case $k<m$, we have found the following
\begin{equation}
t_{m,k}^{(A,+)} = (-1)^{k-m+1} \frac{(2m)!}{2^{3m-k} (m!)^2} \,  \sum_{p=0}^{m-k} \binom{2m+1}{m-k-p} \frac{k^p}{p!}  .
 \end{equation}

\nl 
The case $k=m$ follow easily from the previous discussions. Combining everything together, we can write the asymptotic formula in \ref{ASINTOTICA} for $m \ge 2$ with 
\begin{equation} \label{tmkA}
 t_{m,k}^{(A)} = 
    \begin{cases}
        - \frac{(2m)!}{2^{m+k} (m!)^2} \frac{m \, k^{k-m}}{k\, (k-m)!} \quad &\mbox{if} \ \ k> m \ge 2
        \\[5pt]
        (-1)^{k-m+1} \frac{(2m)!}{2^{3m-k} (m!)^2} \,  \sum\limits_{p=0}^{m-k} \binom{2m+1}{m-k-p} \frac{k^p}{p!} \quad &\mbox{if} \ \ m> k, \, m \ge 2\\[5pt]
        - \frac{(2m)!}{2^{2m} (m!)^2}  \quad &\mbox{if} \ \ m= k, \, m \ge 2
    \end{cases}
\end{equation}

\nl
Now we can study what happens for $m^*=m+2$.

\nl
{\bf Case B:} ${\mathbf{m <2 \ \   \mbox{i.e.} \ \  m^*=m+2.}}$

\nl
Starting from \ref{HansenExpansion} but for $m=0,1$, the leading term will be for $n=m+2$, so that
$$
t_{m,k}^{(B)}= - \frac{(2m+2)! }{2^{2m+3} ((m+1)!)^2 } [X^{m+2,m}_k]
$$
and, as we previously did, we have to look to two different subcases.

\nl
Case B.1: ${\mathbf{m-k<0}}$: here we consider $n=m+2$ and $s=k-m$ so that $t=0$ and \ref{balmino} becomes
\begin{equation}
    [X^{m+2,m}_{k}] = \frac{(-1)^{k-m}}{2^{k-m}}  \binom{2m+3}{0}  \sum_{q=0}^{k-m} \binom{3}{k-m-q} \frac{k^q}{q!} (-1)^q \bigg[ 2 \binom{k -2m-q-4}{0} - \binom{k -2m-q-3}{0}   \bigg].
\end{equation}
Notice that in this case $q$ can take all the integer values between $k-m-3$ and $k-m$, so that
$$
[X^{m+2,m}_{k}] = \frac{(-1)^{k-m}}{2^{k-m}}    \sum_{q=k-m-3}^{k-m} \binom{3}{k-m-q} \frac{k^q}{q!} (-1)^q 
$$
where we use the convention that the product of a sequence on an empty set of indices is $1$. So, using the same notation as case 1, we have obtained for the case $m<k$
\begin{equation}
 t_{m,k}^{(B,-)} =- \frac{(2m+2)! }{2^{2m+3} ((m+1)!)^2 } \frac{(-1)^{k-m}}{2^{k-m}}    \sum_{q=k-m-3}^{k-m} \binom{3}{k-m-q} \frac{k^q}{q!} (-1)^q.
\end{equation}

\nl
Case B.2: ${\mathbf{m-k>0}}$:
Now, setting always $s=k-m$, $n=m+2$, as we did before we have to set $t=m-k$.
\begin{equation}
\begin{split}
    [X^{m+2,m}_{k}] = &(-1)^{k-m} \bigg( \frac{1}{2} \bigg)^{m-k} \bigg\{ \sum_{j=0}^{m-k} \sum_{p=0}^j \binom{2m+3}{j-p} \frac{k^p}{p!} \sum_{q=0}^{k-m+j} \binom{3}{k-m+j-q} \frac{k^q}{q!} (-1)^q
    \\
    &\bigg[ 2 \binom{ -k-p-q-4}{m-k-j} - \binom{-k-p-q-3}{m-k-j}   \bigg]  \bigg\}.
\end{split}
\end{equation}
As in the previous case the sum is restricted to $j=m-k$ and $q=0$:
\begin{equation}
\begin{split}
    [X^{m,m}_{k}] &= (-1)^{k-m} \bigg( \frac{1}{2} \bigg)^{m-k} \bigg\{  \sum_{p=0}^{m-k} \binom{2m+3}{m-k-p} \frac{k^p}{p!} \bigg[ 2 \binom{ -k-p-4}{0} - \binom{-k-p-3}{0}   \bigg]  \bigg\}
    \\
    &=\frac{(-1)^{k-m}}{2^{m-k}} \bigg\{  \sum_{p=0}^{m-k} \binom{2m+3}{m-k-p} \frac{k^p}{p!}  \bigg\}.
\end{split}
\end{equation}
so that for $m=0,1$
\begin{equation}
t_{m,k}^{(B,+)} = (-1)^{k-m+1} \frac{(2m+2)! }{2^{3m-k+3} ((m+1)!)^2 } \bigg\{  \sum_{p=0}^{m-k} \binom{2m+3}{m-k-p} \frac{k^p}{p!}  \bigg\}.
\end{equation}

\nl
In this way we can write the asymptotic coefficient for $m=0,1$ as
\begin{equation} \label{tmkB}
 t_{m,k}^{(B)} = 
    \begin{cases}(-1)^{k-m}
    \frac{(2m+2)!}{2^{m+k+3} ((m+1)!)^2 }    \displaystyle\sum_{q=k-m-3}^{k-m} \binom{3}{k-m-q} \frac{k^q}{q!} (-1)^q \quad &\mbox{if} \ \ k> m, \, m \in \{ 0,1\}
        \\[6pt]
        (-1)^{k-m+1} \frac{(2m+2)! }{2^{3m-k+3} ((m+1)!)^2 } \bigg\{  \displaystyle\sum_{p=0}^{m-k} \binom{2m+3}{m-k-p} \frac{k^p}{p!}  \bigg\} \quad &\mbox{if} \ \ m> k, , \, m \in \{ 0,1\} \\[6pt]
        - \frac{(2m+2)! }{2^{2m+3} ((m+1)!)^2 } \quad &\mbox{if} \ \ m= k, , \, m \in \{ 0,1\}
    \end{cases}
\end{equation}

To be more clear, using the previous notations, we summarize what we have computed with the following expression of the asymptotic for the Fourier coefficient $f_{m,k}$ for small values of $(a,e)$:
\begin{equation}
\label{ASINTOTICA2}
    f_{m,k}(a,e)=
    \begin{cases}
        \mbox{if} \ (m,k)=(0,0) \quad -1+ t_{0,
        0}^{(B)}\,  \,a^{2} \, \big[ 1 + \mathcal{O}(e^2;a) \big] \\[7pt]
        \mbox{if} \ (m,k)\neq(0,0) \quad\begin{cases} \mbox{if} \ m \in \{0,1\} \quad 
    2\,t_{m,k}^{(B)}\, e^{|m-k|} \,a^{m + 2} \, \big[ 1 + \mathcal{O}(e^2;a) \big] \\[5pt]
    \mbox{if} \ m \ge 2 \quad 2\,t_{m,k}^{(A)}\, e^{|m-k|} \,a^{m} \, \big[ 1 + \mathcal{O}(e^2;a) \big] 
    \end{cases} 
    \end{cases}
\end{equation}
then depending on whether k is greater or less than m we will choose the expression of $t_{m,k}^{(\cdot,\pm)}$.

\section{On the Hansen coefficients} \label{HANSEN}

In this section we are going to review the standard theory about Hansen coefficients used in the expansion shown in \equ{HansenExpansion} of the PCR3BP perturbing function.
\smallskip

\nl
Hansen coefficients (Cefola \cite{ce77}) constitute a fundamental class of functions in Celestial Mechanics, appearing notably in planetary theory (Newcomb \cite{ne95}) and in the study of artificial satellite motion (Allan \cite{al67}; Hughes \cite{hu77}). Generalized expansions of these coefficients (Klioner et al.~\cite{kl98}; Sharaf \cite{sh85}, \cite{sh10}, Wu \& Zhang \cite{cinesi2}) have proven essential for describing elliptic motion.

\nl
Giacaglia (\cite{gi76}) was the first to recognize the role of Hansen coefficients in satellite theory, identifying their presence in the disturbing function due to the primary and third bodies. He developed recurrence relations for these coefficients and their derivatives, crucial for evaluating perturbations from geopotential and third-body effects. A similar occurrence of Hansen coefficients has been identified in the disturbing function of 3BP (\cite{LB10}).

\nl
In \cite{gia87}, Giacaglia further demonstrated that Hansen coefficients in the Fourier series with respect to mean anomaly correspond to a rotation of the orbital plane proportional to the eccentricity. These coefficients are expressed through Bessel functions and generalized associated Legendre functions, arising from the rotation of spherical harmonics. Hughes (\cite{hug81}) provided extensive tables of analytical expressions for $X_{0}^{n,\pm m}(e)$ and $X_{0}^{-(n+1),\pm m}(e)$ for $1\leq n\leq 30$ and $0\leq m\leq n$.

\nl
Branham (\cite{bra90}) introduced recursive methods for computing Hansen coefficients via Tisserand’s method and the Von Zeipel-Andoyer method, both in explicit and recursive polynomial forms. Vakhidov (\cite{vak00}) proposed efficient polynomial approximations for these coefficients as functions of eccentricity.

\nl
Hansen coefficients were also employed by He and Zhang (\cite{mia90}) to model perturbations on Flora group asteroids due to Jupiter. Breiter et al. (\cite{bre04}) extended the theory to generalized coefficients $X_{k}^{\gamma j}$ for real $\gamma$, enabling their application to perturbed problems involving drag forces.

\nl
Sadov (\cite{sad08}) deals analytically with the properties of Hansen's coefficients in the theory of elliptic motion considered as functions of the parameter
$\eta=\sqrt{1-e^{2}}$ where $e$ is the eccentricity.

In the next sections, using recursive relations already known in the literature we provide tables with truncated coefficients $X^{n,m}_k$ with also high values of $k$, in particular for $k=0,1,4,8,10.$

\subsection*{Computation of $X_0^{n,m}$ and $X_0^{-(n+1),m}$}

As we have seen before in \ref{Giacaglia} the general Hansen coefficient $X_k^{n,m} (e)$ is a function of the orbital eccentricity and is defined by the generating function
\begin{equation}
        \bigg(\frac{r}{a}\bigg)^n\, \mbox{exp}\,({im f}) = \sum_{k=-\infty}^{+\infty} X^{n,m}_k (e) \ 
        \bigg(\frac{r}{a}\bigg)^n\, \mbox{exp}\,({ik \ell}).
\end{equation}
where $n,m$ and $k$ are integers which may be positive or negative, $r$ the radius vector, $a$ the semi-major axis, $e$ the orbital eccentricity, $f$ the true anomaly and $\ell$ the mean anomaly. The individual coefficients being given by the integral
\begin{equation} \label{IntegraliHansen}
    X_k^{n,m} (e) = \frac{1}{2 \pi} \int_0^{2 \pi} \bigg(\frac{r}{a} \bigg)^n \cos (m f - k \ell) d \ell.
\end{equation}
that shows easily that $X_k^{n,m} = X_{-k}^{n,-m}$.

A number of authors have given extensive table of  these coefficients, the most important are by Cayley (1861, \cite{Cay61}), Newcomb (1895, \cite{ne95}) and Cherniack (1972, \cite{Ch72}) but they are quite tedious and time consuming. We prefer, following the expansion by Hughes \cite{hug81}, to begin with computing $X_0^{n,m}$ and $X_0^{-(n+1),m}$ for $0\le n \le 15$ and $0 \le m \le 3$ (see table \ref{res0}).
\smallskip

\nl
If we put $k=0$, then the integrals \ref{IntegraliHansen} for $X_0^{n,m}$ and $X_0^{-(n+1),m}$ become
\begin{equation} \label{Hughes}
    X_0^{n,m} = \frac{1}{2\pi} \int_0^{2\pi} \bigg(\frac{r}{a} \bigg)^n \cos (m f) d \ell, \qquad X_0^{-(n+1),m} = \frac{1}{2\pi} \int_0^{2\pi} \bigg(\frac{a}{r} \bigg)^{n+1} \cos (m f) d \ell.
\end{equation}
On putting $m=-m$ into \ref{Hughes} it is obvious that $X_0^{n,m} = X_0^{n,-m}$ and $X_0^{-(n+1),m} = X_0^{-(n+1),-m}$, therefore it is only necessary to obtain relations for positive $m$. 
If the integrals \ref{Hughes} are evalueted (see for example \cite{hu77}) we have
\begin{equation}
    \begin{split}
        X_0^{n,m} &= \bigg(-\frac{e}{2} \bigg)^m \binom{n+m+1}{m} F \bigg( \frac{m-n-1}{2}, \frac{m-n}{2}, m+1 ;e^2   \bigg),
        \\
        X_0^{-(n+1),m} &= \bigg(-\frac{e}{2} \bigg)^m \frac{1}{(1-e)^{(2n-1)/2}} \sum_{j=0}^{[(n-m-1)/2]} \frac{1}{2^j} \binom{n-1}{2j+m} \binom{2j+m}{j} e^2j,
    \end{split}
\end{equation}
where $F( \ \ )$ is the standard hypergeometric function and $[ \ \ ]$ denotes the nearest lowest integer.
From these equations, replacing the hypergeometric functional expressions, one can obtain the recursive formulae (see \cite{br75}) for $X_0^{n,m}$
\begin{equation}
    \begin{split}
    X_0^{n+1,m} = \frac{(2n+3)}{(n+2)} X_0^{n,m} - \frac{(n+1-m)(n+1+m)}{(n+1)(n+2)} (1-e^2) X_0^{n-1,m}
    \\
    e X_0^{n,m+1} = \frac{1}{(n-m+1)} (e(n+m+1) X_0^{n,m-1} + 2m X_0^{n,m}),
    \end{split}
\end{equation}
the corresponding recoursive relations for $X_0^{-(n+1),m}$ are
\begin{equation}
    \begin{split}
        (n-m+1)(n-m-1) &X_0^{-(n+3),m} = \frac{(n+1)}{(1-e^2)} [(2n+1)X_0^{-(n+2),m} - n X_0^{-(n+1),m}]
        \\
        X_0^{-(n+1),m} &= \frac{1}{(n-m-1)} [ 2(m+1)\sqrt{1-e^2} X_0^{-(n+1),(m+1)} +
        \\
        & + (n+m+1)e^2 (1-e^2)^{3/2} X_0^{-(n+1),(m+2)}].
    \end{split}
\end{equation}


\begin{table}[htbp]
\centering 
\arraycolsep=6pt
\renewcommand{\arraystretch}{2}

\resizebox{\textwidth}{!}{
  $\displaystyle          
\begin{array}{|r|l|l|l|l|}
\toprule
     n & X^{n,0}_{0} 
      & X^{n,1}_{0}
      & X^{n,2}_{0}
      & X^{n,3}_{0}\\
\midrule
0 &
  1 &
  -e &
  \dfrac{3e^{2}}{4}+\dfrac{e^{4}}{8}+\dfrac{3e^{6}}{64} &
  -\dfrac{e^{3}}{2}-\dfrac{3e^{5}}{16}-\dfrac{3e^{7}}{32}\\[2pt]

1 &
  1+\dfrac{e^{2}}{2} &
  -\dfrac{3e}{2} &
  \dfrac{3e^{2}}{2} &
  -\dfrac{5e^{3}}{4}-\dfrac{5e^{5}}{32}-\dfrac{3e^{7}}{64}\\

2 &
  1+\dfrac{3e^{2}}{2} &
  -2e-\dfrac{e^{3}}{2} &
  \dfrac{5e^{2}}{2} &
  -\dfrac{5e^{3}}{2}\\

3 &
  1+3e^{2}+\dfrac{3e^{4}}{8} &
  -\dfrac{5e}{2}-\dfrac{15e^{3}}{8} &
  \dfrac{15e^{2}}{4}+\dfrac{5e^{4}}{8} &
  -\dfrac{35e^{3}}{8}\\

4 &
  1+5e^{2}+\dfrac{15e^{4}}{8} &
  -3e-\dfrac{9e^{3}}{2}-\dfrac{3e^{5}}{8} &
  \dfrac{21e^{2}}{4}+\dfrac{21e^{4}}{8} &
  -7e^{3}-\dfrac{7e^{5}}{8}\\

5 &
  1+\dfrac{15e^{2}}{2}+\dfrac{45e^{4}}{8}+\dfrac{5e^{6}}{16} &
  -\dfrac{7e}{2}-\dfrac{35e^{3}}{4}-\dfrac{35e^{5}}{16} &
  7e^{2}+7e^{4}+\dfrac{7e^{6}}{16} &
  -\dfrac{21e^{3}}{2}-\dfrac{63e^{5}}{16}\\

6 &
  1+\dfrac{21e^{2}}{2}+\dfrac{105e^{4}}{8}+\dfrac{35e^{6}}{16} &
  -4e-15e^{3}-\dfrac{15e^{5}}{2}-\dfrac{5e^{7}}{16} &
  9e^{2}+15e^{4}+\dfrac{45e^{6}}{16} &
  -15e^{3}-\dfrac{45e^{5}}{4}-\dfrac{9e^{7}}{16}\\

7 &
  1+14e^{2}+\dfrac{105e^{4}}{4}+\dfrac{35e^{6}}{4} &
  -\dfrac{9e}{2}-\dfrac{189e^{3}}{8}-\dfrac{315e^{5}}{16}-\dfrac{315e^{7}}{128} &
  \dfrac{45e^{2}}{4}+\dfrac{225e^{4}}{8}+\dfrac{675e^{6}}{64} &
  -\dfrac{165e^{3}}{8}-\dfrac{825e^{5}}{32}-\dfrac{495e^{7}}{128}\\

8 &
  1+18e^{2}+\dfrac{189e^{4}}{4}+\dfrac{105e^{6}}{4} &
  -5e-35e^{3}-\dfrac{175e^{5}}{4}-\dfrac{175e^{7}}{16} &
  \dfrac{55e^{2}}{4}+\dfrac{385e^{4}}{8}+\dfrac{1925e^{6}}{64} &
  -\dfrac{55e^{3}}{2}-\dfrac{825e^{5}}{16}-\dfrac{495e^{7}}{32}\\

9 &
  1+\dfrac{45e^{2}}{2}+\dfrac{315e^{4}}{4}+\dfrac{525e^{6}}{8} &
  -\dfrac{11e}{2}-\dfrac{99e^{3}}{2}-\dfrac{693e^{5}}{8}-\dfrac{1155e^{7}}{16} &
  \dfrac{33e^{2}}{2}+77e^{4}+\dfrac{1155e^{6}}{16} &
  -\dfrac{143e^{3}}{4}-\dfrac{3003e^{5}}{32}-\dfrac{3003e^{7}}{64}\\

10 &
  1+\dfrac{55e^{2}}{2}+\dfrac{495e^{4}}{4}+\dfrac{1155e^{6}}{8} &
  -6e-\dfrac{135e^{3}}{2}-\dfrac{315e^{5}}{2}-\dfrac{1575e^{7}}{16} &
  \dfrac{39e^{2}}{2}+117e^{4}+\dfrac{2457e^{6}}{16} &
  -\dfrac{91e^{3}}{2}-\dfrac{637e^{5}}{4}-\dfrac{1911e^{7}}{16}\\

11 &
  1+{33e^{2}}+\dfrac{1485e^{4}}{8}+\dfrac{1155e^{6}}{4} &
  -\dfrac{13e}{2}-\dfrac{715e^{3}}{8}-\dfrac{2145e^{5}}{8}-\dfrac{15015e^{7}}{64} &
  \dfrac{91e^{2}}{4}+\dfrac{1365e^{4}}{8}+\dfrac{9555e^{6}}{32} &
  -\dfrac{455e^{3}}{8}-\dfrac{4095e^{5}}{16}-\dfrac{17199e^{7}}{64}\\

12 &
  1+{39e^{2}}+\dfrac{2145e^{4}}{8}+\dfrac{2145e^{6}}{4} &
  -7e-\dfrac{231e^{3}}{2}-\dfrac{3465e^{5}}{8}-\dfrac{8085e^{7}}{16} &
  \dfrac{105e^{2}}{4}+\dfrac{1925e^{4}}{8}+\dfrac{17325e^{6}}{32} &
  -{70e^{3}}-\dfrac{1575e^{5}}{4}-\dfrac{2205e^{7}}{4}\\

13 &
  1+\dfrac{91e^{2}}{2}+\dfrac{3003e^{4}}{8}+\dfrac{15015e^{6}}{16} &
  -\dfrac{15e}{2}-\dfrac{585e^{3}}{4}-\dfrac{10725e^{5}}{16}-\dfrac{32175e^{7}}{32} &
 30e^{2}+330e^{4}+\dfrac{7425e^{6}}{8} &
  -85e^{3}-\dfrac{4675e^{5}}{8}-\dfrac{8415e^{7}}{8}\\

14 &
  1+\dfrac{105e^{2}}{2}+\dfrac{4095e^{4}}{8}+\dfrac{25025e^{6}}{16} &
  -8e-182e^{3}-1001e^{5}-\dfrac{15015e^{7}}{8} &
  34e^{2}+442e^{4}+\dfrac{12155e^{6}}{8} &
  -102e^{3}-\dfrac{1683e^{5}}{2}-\dfrac{15147e^{7}}{8}\\

15 &
  1+60e^{2}+\dfrac{1365e^{4}}{2}+\dfrac{5005e^{6}}{2} &
  -\dfrac{17e}{2}-\dfrac{1785e^{3}}{8}-\dfrac{23205e^{5}}{16}-\dfrac{425425e^{7}}{128} &
  \dfrac{153e^{2}}{4}+\dfrac{4641e^{4}}{8}+\dfrac{153153e^{6}}{64} &
  -\dfrac{969e^{3}}{8}-\dfrac{37791e^{5}}{32}-\dfrac{415701e^{7}}{128}\\
\bottomrule
  \end{array}
  $
}
\caption{Hansen coefficients $X_0^{n,m} $ for  $0 \leq n \leq 15$, $m=0,1,2,3$  at order $7$ in $e$.}
\label{res0}
\end{table}
\FloatBarrier

\subsection*{Computation of $X_k^{n,m}$ and $X_k^{-(n+1),m}$ when $k \ne 0$.}
If $k \ne 0$ then the computation of $X_k^{n,m}$ and $X_k^{-(n+1),m}$ presents some difficulty in that the analytical expressions for such coefficients do not terminate, consequently the series have to be truncated at some particular order in the eccentricity.

\nl
Since most planets and satellites both natural and artificial have small or moderate eccentricities ($0 \le e \le 0.1$), a series expansion in the eccentricity is usually fine.
\begin{equation}\label{espansione1}
    X_k^{n,m} = \sum_{q} \widehat{X}_{k,q}^{n,m} e^{q}
\end{equation}
The coefficients $\widehat{X}_{k,q}^{n,m}$ with shifted indices are known as Newcomb’s operators defined as follows
\begin{equation}
    X_{\rho,\sigma}^{n,m} = \widehat{X}_{m+\rho-\sigma,\rho+\sigma}^{n,m}
\end{equation}
in such a way that the expansion in \ref{espansione1} becomes  the well-known
\begin{equation}\label{espansione hansen}
    X_k^{n,m} = \sum_{\rho-\sigma = k +m} X_{\rho,\sigma}^{n,m} e^{\rho + \sigma}.
\end{equation}
For $\sigma=0$, knowing that $X^{n,m}_{0,0}=1$ and $X^{n,m}_{1,0}=\bigg( m-\frac{n}{2} \bigg)$, the recursive relations are easily founded 
\begin{equation}
    4 \rho X^{n,m}_{\rho,0} = 2(2m-n) X_{\rho-1,0}^{n,m+1} + (m-n) X_{\rho-2,0}^{n,m+2},
\end{equation}
while for $\sigma \ne 0$ the relation
\begin{equation}
    \begin{split}
        4 \sigma X_{\rho,\sigma}^{n,m} = &- 2(2m+n) X_{\rho,\sigma-1}^{n,m-1} - (m+n) X_{\rho,\sigma-2}^{n,m-2} - (\rho - 5\sigma + 4 + 4m + n) X_{\rho-1,\sigma-1}^{n,m} + 
        \\
        &+2(\rho-\sigma+m) \sum_{j \ge 2} (-1)^j \binom{3/2}{j} X_{\rho-j,\sigma-j}^{n,m}
    \end{split}
\end{equation}
is used. From the above relations we notice that $X^{n,m}_{\rho,\sigma}=0$ whenever $\rho$ or $\sigma$ is negative, and that
\begin{equation}
    X_{\rho,\sigma}^{n,m} = X_{\sigma,\rho}^{n,-m} \ \ \ \mbox{if} \ \sigma > \rho.
\end{equation}

\nl
As we can see, when $k \ne0$ the recursive relations are more complicate than the ones for $k=0$. In fact, for computational purposes, the fastest way to calculate these Hansen coefficients with $k \ne 0$ is by using the Bessel function, as shown in \cite{cinesi}.
In particular, we will use the following expansion, known as {\it Wnuk's method}:
\begin{equation}
    \begin{split}
    X_k^{n,m} =& (1+ \beta^2)^{-(n+1)} \sum_{t=-\infty}^\infty E_{k-t}^{n,m} J_t (ke),
    \\
    2(1-e^2) \frac{d X_k^{n,m}}{ d e} =& - \frac{2m}{e} X_k^{n,m} - (n+m)e X_k^{n,m} + \frac{2k (1-e^2)^{3/2}}{e} X_k^{n,m} - 
    \\
    &(2n + 4m) X_k^{n,m-1} - (n+m)e X_k^{n,m-2},
    \end{split}
\end{equation}
where
\begin{equation}
    \begin{split}
        \beta &= \frac{e}{1 + \sqrt{1- e^2}},
        \\
        E_{k-t}^{n,m} &= \begin{cases}
            (-\beta)^{k-t-m} \sum\limits_{s=0}^\infty \binom{n-m+1}{k-t-m+s} \binom{n+m+1}{s} \beta^{2s}, \quad (k-t-m \ge 0), \\
            (-\beta)^{t-k+m} \sum\limits_{s=0}^\infty \binom{n+m+1}{t-k+m+s} \binom{n-m+1}{s} \beta^{2s}, \quad (k-t-m <0 ),
        \end{cases}
    \end{split}
\end{equation}
and $J_t (ke)$ is the Bessel function of $ke$, for which the following relation holds (for $t<0$ or $ke<0$)
$$
J_{-t} (ke) = J_{t} (-ke) = (-1)^t J_{t} (ke).
$$


\begin{table}[htbp]
\centering
\arraycolsep=6pt
\renewcommand{\arraystretch}{2}

\resizebox{\textwidth}{!}{
  $\displaystyle           
\begin{array}{|r|l|l|l|l|}
\toprule
     n & X^{n,0}_{1} 
      & X^{n,1}_{1}
      & X^{n,2}_{1}
      & X^{n,3}_{1}\\
\midrule
0 &
  0 &
  1 - e^{2} + \dfrac{7e^{4}}{64} - \dfrac{5e^{6}}{288} &
  -2e + \dfrac{7e^{3}}{4} - \dfrac{5e^{5}}{96} + \dfrac{271e^{7}}{4608} &
  \dfrac{21e^{2}}{8} - \dfrac{31e^{4}}{16} - \dfrac{103e^{6}}{1024}\\
1 &
  -\dfrac{e}{2} + \dfrac{3e^{3}}{16} - \dfrac{5e^{5}}{384} + \dfrac{7e^{7}}{18432} &
  1 - \dfrac{e^{2}}{2} - \dfrac{e^{4}}{64} - \dfrac{29e^{6}}{1152} &
  -\dfrac{5e}{2} + \dfrac{33e^{3}}{16} - \dfrac{73e^{5}}{384} + \dfrac{881e^{7}}{18432} &
  \dfrac{31e^{2}}{8} - \dfrac{77e^{4}}{24} + \dfrac{155e^{6}}{1024}\\
2 &
  -e + \dfrac{e^{3}}{8} - \dfrac{e^{5}}{192} + \dfrac{e^{7}}{9216} &
  1 + \dfrac{e^{2}}{2} - \dfrac{25e^{4}}{64} - \dfrac{23e^{6}}{1152} &
  -3e + \dfrac{13e^{3}}{8} + \dfrac{5e^{5}}{192} + \dfrac{227e^{7}}{3072} &
  \dfrac{43e^{2}}{8} - \dfrac{25e^{4}}{6} + \dfrac{1069e^{6}}{3072}\\ 
3 &
  -\dfrac{3e}{2} - \dfrac{9e^{3}}{16} + \dfrac{15e^{5}}{128} - \dfrac{35e^{7}}{6144} &
  1 + 2e^{2} - \dfrac{41e^{4}}{64} - \dfrac{37e^{6}}{576} &
  -\dfrac{7e}{2} + \dfrac{e^{3}}{16} + \dfrac{289e^{5}}{384} + \dfrac{1645e^{7}}{18432} &
  \dfrac{57e^{2}}{8} - \dfrac{65e^{4}}{16} - \dfrac{19e^{6}}{1024}\\ 
4 &
  -2e - \dfrac{9e^{3}}{4} + \dfrac{19e^{5}}{96} - \dfrac{29e^{7}}{4608} &
  1 + 4e^{2} - \dfrac{e^{4}}{64} - \dfrac{199e^{6}}{576} &
  -4e - 3e^{3} + \dfrac{79e^{5}}{48} + \dfrac{233e^{7}}{1152} &
  \dfrac{73e^{2}}{8} - \dfrac{91e^{4}}{48} - \dfrac{1419e^{6}}{1024}\\ 
5 &
  -\dfrac{5e}{2} - \dfrac{85e^{3}}{16} - \dfrac{185e^{5}}{384} + \dfrac{1535e^{7}}{18432} &
  1 + \dfrac{13e^{2}}{2} + \dfrac{167e^{4}}{64} - \dfrac{995e^{6}}{1152} &
  -\dfrac{9e}{2} - \dfrac{127e^{3}}{16} + \dfrac{595e^{5}}{384} + \dfrac{4763e^{7}}{6144} &
  \dfrac{91e^{2}}{8} - \dfrac{43e^{4}}{12} - \dfrac{10979e^{6}}{3072}\\ 
6 &
  -3e - \dfrac{81e^{3}}{8} - \dfrac{225e^{5}}{64} + \dfrac{757e^{7}}{3072} &
  1 + \dfrac{19e^{2}}{2} + \dfrac{559e^{4}}{64} - \dfrac{929e^{6}}{1152} &
  -5e - \dfrac{121e^{3}}{8} - \dfrac{349e^{5}}{192} + \dfrac{19643e^{7}}{9216} &
  \dfrac{111e^{2}}{8} + \dfrac{111e^{4}}{8} - \dfrac{5125e^{6}}{1024}\\ 
7 &
  -\dfrac{7e}{2} - \dfrac{273e^{3}}{16} - \dfrac{4487e^{5}}{384} - \dfrac{5873e^{7}}{18432} &
  1 + 13e^{2} + \dfrac{1295e^{4}}{64} + \dfrac{43e^{6}}{18} &
  -\dfrac{11e}{2} - \dfrac{399e^{3}}{16} - \dfrac{4675e^{5}}{384} + \dfrac{66317e^{7}}{18432} &
  \dfrac{133e^{2}}{8} + \dfrac{1475e^{4}}{48} - \dfrac{1727e^{6}}{1024}\\ 
8 &
  -4e - \dfrac{53e^{3}}{2} - \dfrac{1405e^{5}}{48} - \dfrac{10745e^{7}}{2304} &
  1 + 17e^{2} + \dfrac{2519e^{4}}{64} + \dfrac{2057e^{6}}{144} &
  -6e - \dfrac{151e^{3}}{4} - \dfrac{3359e^{5}}{96} + \dfrac{825e^{7}}{512} &
  \dfrac{157e^{2}}{8} + \dfrac{2695e^{4}}{48} + \dfrac{43123e^{6}}{3072}\\ 
9 &
  -\dfrac{9e}{2} - \dfrac{621e^{3}}{16} - \dfrac{7983e^{5}}{128} - \dfrac{42209e^{7}}{2048} &
  1 + \dfrac{43e^{2}}{2} + \dfrac{4399e^{4}}{64} + \dfrac{51847e^{6}}{1152} &
  -\dfrac{13e}{2} - \dfrac{863e^{3}}{16} - \dfrac{29873e^{5}}{384} - \dfrac{246767e^{7}}{18432} &
  \dfrac{183e^{2}}{8} + \dfrac{739e^{4}}{8} + \dfrac{56291e^{6}}{1024}\\ 
10 &
  -5e - \dfrac{435e^{3}}{8} - \dfrac{22885e^{5}}{192} - \dfrac{595315e^{7}}{9216} &
  1 + \dfrac{53e^{2}}{2} + \dfrac{7127e^{4}}{64} + \dfrac{128005e^{6}}{1152} &
  -7e - \dfrac{591e^{3}}{8} - \dfrac{28895e^{5}}{192} - \dfrac{567203e^{7}}{9216} &
  \dfrac{211e^{2}}{8} + \dfrac{1703e^{4}}{12} + \dfrac{144167e^{6}}{1024}\\ 
11 &
  -\dfrac{11e}{2} - \dfrac{1177e^{3}}{16} - \dfrac{80795e^{5}}{384} - \dfrac{3080473e^{7}}{18432} &
  1 + 32e^{2} + \dfrac{10919e^{4}}{64} + \dfrac{136973e^{6}}{576} &
  -\dfrac{15e}{2} - \dfrac{1567e^{3}}{16} - \dfrac{102023e^{5}}{384} - \dfrac{1105553e^{7}}{6144} &
  \dfrac{241e^{2}}{8} + \dfrac{9961e^{4}}{48} + \dfrac{921535e^{6}}{3072}\\ 
12 &
  -6e - \dfrac{387e^{3}}{4} - \dfrac{11181e^{5}}{32} - \dfrac{583829e^{7}}{1536} &
  1 + 38e^{2} + \dfrac{16015e^{4}}{64} + \dfrac{265751e^{6}}{576} &
  -8e - \dfrac{253e^{3}}{2} - \dfrac{10535e^{5}}{24} - \dfrac{992941e^{7}}{2304} &
  \dfrac{273e^{2}}{8} + \dfrac{4675e^{4}}{16} + \dfrac{585773e^{6}}{1024}\\ 
13 &
  -\dfrac{13e}{2} - \dfrac{1989e^{3}}{16} - \dfrac{212225e^{5}}{384} - \dfrac{14492465e^{7}}{18432} &
  1 + \dfrac{89e^{2}}{2} + \dfrac{22679e^{4}}{64} + \dfrac{958561e^{6}}{1152} &
  -\dfrac{17e}{2} - \dfrac{2559e^{3}}{16} - \dfrac{264661e^{5}}{384} - \dfrac{16828495e^{7}}{18432} &
  \dfrac{307e^{2}}{8} + \dfrac{2395e^{4}}{6} + \dfrac{1034151e^{6}}{1024}\\ 
14 &
  -7e - \dfrac{1253e^{3}}{8} - \dfrac{161287e^{5}}{192} - \dfrac{13926773e^{7}}{9216} &
  1 + \dfrac{103e^{2}}{2} + \dfrac{31199e^{4}}{64} + \dfrac{1632667e^{6}}{1152} &
  -9e - \dfrac{1589e^{3}}{8} - \dfrac{199489e^{5}}{192} - \dfrac{5448827e^{7}}{3072} &
  \dfrac{343e^{2}}{8} + \dfrac{12767e^{4}}{24} + \dfrac{5168953e^{6}}{3072}\\ 
15 &
  -\dfrac{15e}{2} - \dfrac{3105e^{3}}{16} - \dfrac{158085e^{5}}{128} - \dfrac{16819915e^{7}}{6144} &
  1 + 59e^{2} + \dfrac{41887e^{4}}{64} + \dfrac{663745e^{6}}{288} &
  -\dfrac{19e}{2} - \dfrac{3887e^{3}}{16} - \dfrac{581675e^{5}}{384} - \dfrac{59463251e^{7}}{18432} &
  \dfrac{381e^{2}}{8} + \dfrac{11109e^{4}}{16} + \dfrac{2742377e^{6}}{1024}\\
\bottomrule
  \end{array}
  $
}
\caption{Hansen coefficients $X_1^{n,m} $ for  $0 \leq n \leq 15$, $m=0,1,2,3$ at order $7$ in $e$.}
\label{res1}
\end{table}
\FloatBarrier

\begin{table}[htbp]
\centering
\arraycolsep=6pt
\renewcommand{\arraystretch}{2}

\resizebox{\textwidth}{!}{
  $\displaystyle           
\begin{array}{|r|l|l|l|l|}
\toprule
     n & X^{n,0}_{4} 
      & X^{n,1}_{4}
      & X^{n,2}_{4}
      & X^{n,3}_{4}\\
\midrule
0 & 0 
  & \displaystyle \frac{4e^3}{3} - \frac{7e^5}{3} + \frac{83e^7}{60}
  & \displaystyle \frac{13e^2}{4} - \frac{259e^4}{24} + \frac{559e^6}{48}
  & \displaystyle 3e - \frac{39e^3}{2} + \frac{155e^5}{4} - \frac{767e^7}{24} \\

1 & \displaystyle -\frac{e^4}{6} + \frac{e^6}{5}
  & \displaystyle \frac{e^3}{3} - \frac{5e^5}{12} + \frac{43e^7}{240}
  & \displaystyle 2e^2 - \frac{19e^4}{3} + \frac{55e^6}{8}
  & \displaystyle \frac{5e}{2} - \frac{63e^3}{4} + \frac{757e^5}{24} - \frac{1979e^7}{72} \\

2 & \displaystyle -\frac{e^4}{12} + \frac{e^6}{15}
  & \displaystyle -\frac{e^3}{6} + \frac{13e^5}{24} - \frac{47e^7}{96}
  & \displaystyle e^2 - \frac{5e^4}{2} + \frac{101e^6}{48}
  & \displaystyle 2e - \frac{23e^3}{2} + \frac{65e^5}{3} - \frac{5317e^7}{288} \\

3 & \displaystyle \frac{e^4}{16} - \frac{3e^6}{20}
  & \displaystyle -\frac{7e^3}{24} + \frac{55e^5}{96} - \frac{701e^7}{1920}
  & \displaystyle \frac{e^2}{4} + \frac{5e^4}{24} - \frac{227e^6}{192}
  & \displaystyle \frac{3e}{2} - \frac{57e^3}{8} + 11e^5 - \frac{2905e^7}{384} \\

4 & \displaystyle \frac{7e^4}{48} - \frac{e^6}{5}
  & \displaystyle -\frac{e^3}{6} + \frac{e^5}{48} + \frac{119e^7}{480}
  & \displaystyle -\frac{e^2}{4} + \frac{37e^4}{24} - \frac{139e^6}{64}
  & \displaystyle e - 3e^3 + \frac{37e^5}{24} + \frac{251e^7}{144} \\

5 & \displaystyle \frac{5e^4}{48} + \frac{e^6}{96}
  & \displaystyle \frac{e^3}{12} - \frac{59e^5}{96} + \frac{673e^7}{960}
  & \displaystyle -\frac{e^2}{2} + \frac{3e^4}{2} - \frac{23e^6}{24}
  & \displaystyle \frac{e}{2} + \frac{e^3}{2} - \frac{235e^5}{48} + \frac{1841e^7}{288} \\

6 & \displaystyle -\frac{e^4}{16} + \frac{57e^6}{160}
  & \displaystyle -\frac{e^3}{3} + \frac{5e^5}{6} + \frac{49e^7}{120}
  & \displaystyle -\frac{e^2}{2} + \frac{e^4}{3} + \frac{35e^6}{24}
  & \displaystyle 3e^3 - 7e^5 + \frac{39e^7}{8} \\

7 & \displaystyle -\frac{7e^4}{24} + \frac{21e^6}{40}
  & \displaystyle \frac{11e^3}{24} + \frac{7e^5}{24} - \frac{139e^7}{192}
  & \displaystyle -\frac{e^2}{4} - \frac{35e^4}{24} + \frac{109e^6}{32}
  & \displaystyle -\frac{e}{2} + \frac{33e^3}{8} - \frac{205e^5}{48} - \frac{1069e^7}{576} \\

8 & \displaystyle -\frac{11e^4}{24} + \frac{17e^6}{120}
  & \displaystyle \frac{e^3}{3} + \frac{25e^5}{24} - \frac{241e^7}{120}
  & \displaystyle \frac{e^2}{4} - \frac{25e^4}{8} + \frac{289e^6}{96}
  & \displaystyle -e + \frac{7e^3}{2} + \frac{8e^5}{3} - \frac{1445e^7}{144} \\
\bottomrule
  \end{array}
  $
}
\caption{Hansen coefficients $X_4^{n,m} $ for  $0 \leq n \leq 8$, $m=0,1,2,3$ at order $7$ in $e$.}
\label{res4}
\end{table}
\FloatBarrier

\begin{table}[!htbp]
\centering
\arraycolsep=6pt
\renewcommand{\arraystretch}{2}

\resizebox{\textwidth}{!}{
  $\displaystyle           
\begin{array}{|r|l|l|l|l|}
\toprule
     n & X^{n,0}_{8} 
      & X^{n,1}_{8}
      & X^{n,2}_{8}
      & X^{n,3}_{8}\\
\midrule
0 & 0 
  & \displaystyle \frac{1024e^7}{315} - \frac{2816e^9}{315}
  & \displaystyle \frac{42037e^6}{2880} - \frac{2321957e^8}{40320} + \frac{2116733e^{10}}{23040}
  & \displaystyle \frac{2611e^5}{80} - \frac{87599e^7}{480} + \frac{155981e^9}{384} \\

1 & \displaystyle -\frac{64e^8}{315} + \frac{256e^{10}}{567}
  & \displaystyle \frac{128e^7}{315} - \frac{32e^9}{35}
  & \displaystyle \frac{256e^6}{45} - \frac{1408e^8}{63} + \frac{11408e^{10}}{315}
  & \displaystyle \frac{8551e^5}{480} - \frac{288221e^7}{2880} + \frac{6112597e^9}{26880} \\

2 & \displaystyle -\frac{16e^8}{315} + \frac{256e^{10}}{2835}
  & \displaystyle -\frac{32e^7}{105} + \frac{344e^9}{315}
  & \displaystyle \frac{64e^6}{45} - \frac{1504e^8}{315} + \frac{2084e^{10}}{315}
  & \displaystyle \frac{128e^5}{15} - \frac{416e^7}{9} + \frac{32384e^9}{315} \\

3 & \displaystyle \frac{2e^8}{35} - \frac{32e^{10}}{189}
  & \displaystyle -\frac{68e^7}{315} + \frac{193e^9}{315}
  & \displaystyle -\frac{8e^6}{45} + \frac{68e^8}{45} - \frac{2329e^{10}}{630}
  & \displaystyle \frac{16e^5}{5} - \frac{76e^7}{5} + \frac{3104e^9}{105} \\

4 & \displaystyle \frac{17e^8}{315} - \frac{352e^{10}}{2835}
  & \displaystyle \frac{2e^7}{315} - \frac{13e^9}{70}
  & \displaystyle -\frac{4e^6}{9} + \frac{622e^8}{315} - \frac{4357e^{10}}{1260}
  & \displaystyle \frac{8e^5}{15} - \frac{26e^7}{45} - \frac{116e^9}{35} \\

5 & \displaystyle -\frac{e^8}{504} + \frac{29e^{10}}{567}
  & \displaystyle \frac{47e^7}{420} - \frac{2221e^9}{5040}
  & \displaystyle -\frac{19e^6}{90} + \frac{629e^8}{1260} - \frac{1019e^{10}}{10080}
  & \displaystyle -\frac{7e^5}{15} + \frac{671e^7}{180} - \frac{92e^9}{9} \\

6 & -\displaystyle \frac{47e^8}{1120} + \frac{26e^{10}}{189}
  & \displaystyle \frac{403e^7}{5040} - \frac{3467e^9}{20160}
  & \displaystyle \frac{23e^6}{360} - \frac{725e^8}{1008} + \frac{78959e^{10}}{40320}
  & \displaystyle -\frac{11e^5}{20} + \frac{667e^7}{240} - \frac{2077e^9}{420} \\

7 & \displaystyle -\frac{403e^8}{11520} + \frac{383e^{10}}{6480}
  & \displaystyle -\frac{563e^7}{40320} + \frac{3947e^9}{17920}
  & \displaystyle \frac{541e^6}{2880} - \frac{36667e^8}{40320} + \frac{69443e^{10}}{46080}
  & \displaystyle -\frac{119e^5}{480} + \frac{307e^7}{1152} + \frac{509e^9}{240} \\

8 & \displaystyle \frac{563e^8}{80640} - \frac{437e^{10}}{4536}
  & \displaystyle -\frac{47e^7}{560} - \frac{29333e^9}{80640}
  & \displaystyle \frac{403e^6}{2880} + \frac{10727e^8}{40320} - \frac{127873e^{10}}{322560}
  & \displaystyle \frac{23e^5}{240} - \frac{2273e^7}{1440} + \frac{212161e^9}{40320} \\
\bottomrule
  \end{array}
  $
}
\caption{Hansen coefficients $X_8^{n,m} $ for  $0 \leq n \leq 8$, $m=0,1,2,3$ at order $10$ in $e$.}
\label{res8}
\end{table}
\FloatBarrier

\begin{table}[!htbp]
\centering
\arraycolsep=6pt
\renewcommand{\arraystretch}{2}

\resizebox{\textwidth}{!}{
  $\displaystyle           
\begin{array}{|r|l|l|l|l|}
\toprule
     n & X^{n,0}_{10} 
      & X^{n,1}_{10}
      & X^{n,2}_{10}
      & X^{n,3}_{10}\\
\midrule
0 & 0 
  & \displaystyle \frac{390625   e^9}{72576} - \frac{5078125   e^{11}}{290304}
  & \displaystyle \frac{461843   e^8}{16128} - \frac{1136803   e^{10}}{9072} + \frac{539506157   e^{12}}{2322432}
  & \displaystyle \frac{106469   e^7}{1344} - \frac{1242377   e^9}{2688} + \frac{218081281   e^{11}}{193536} \\

1 & \displaystyle -\frac{78125   e^{10}}{290304} + \frac{390625   e^{12}}{532224} 
  & \displaystyle \frac{78125   e^9}{145152} - \frac{859375   e^{11}}{580608}
  & \displaystyle \frac{78125   e^8}{8064} - \frac{6171875   e^{10}}{145152} + \frac{31015625   e^{12}}{387072}
  & \displaystyle \frac{305593   e^7}{8064} - \frac{1200935   e^9}{5376} + \frac{646152725   e^{11}}{1161216} \\

2 & \displaystyle -\frac{15625   e^{10}}{290304} + \frac{390625   e^{12}}{3193344}
  & \displaystyle -\frac{15625   e^9}{36288} + \frac{15625   e^{11}}{9072} 
  & \displaystyle \frac{15625   e^8}{8064} - \frac{359375   e^{10}}{48384} + \frac{14234375   e^{12}}{1161216} 
  & \displaystyle \frac{15625   e^7}{1008} - \frac{359375   e^9}{4032} + \frac{31796875   e^{11}}{145152} \\

3 & \displaystyle \frac{3125   e^{10}}{48384} - \frac{78125   e^{12}}{354816}
  & \displaystyle -\frac{34375   e^9}{145152} + \frac{228125   e^{11}}{290304}
  & \displaystyle -\frac{3125   e^8}{8064} + \frac{96875   e^{10}}{36288} - \frac{7878125   e^{12}}{1161216}
  & \displaystyle \frac{3125   e^7}{672} - \frac{128125   e^9}{5376} + \frac{5046875   e^{11}}{96768}
\\

4 & \displaystyle \frac{6875   e^{10}}{145152} - \frac{15625   e^{12}}{118272} 
  & \displaystyle \frac{625   e^9}{18144} - \frac{41875   e^{11}}{145152}
  & \displaystyle -\frac{625   e^8}{1152} + \frac{186875   e^{10}}{72576} - \frac{1980625   e^{12}}{387072}
  & \displaystyle \frac{625  e^7}{2016} + \frac{625  e^9}{672} - \frac{2434375  e^{11}}{290304} \\

5 & \displaystyle -\frac{625   e^{10}}{72576} + \frac{203125   e^{12}}{3193344}
  & \displaystyle \frac{8125   e^9}{72576} - \frac{34375   e^{11}}{72576}
  & \displaystyle -\frac{625  e^8}{4032} + \frac{625  e^{10}}{2016} + \frac{78125  e^{12}}{290304}
  & \displaystyle -\frac{3125   e^7}{4032} + \frac{44375   e^9}{8064} - \frac{8969375   e^{11}}{580608} \\

6 & -\displaystyle \frac{1625   e^{10}}{48384} + \frac{14375   e^{12}}{118272}
  & \displaystyle \frac{3625   e^9}{72576} - \frac{29875  e^{11}}{290304}
  & \displaystyle \frac{125   e^8}{1008} - \frac{67625   e^{10}}{72576} + \frac{1453625   e^{12}}{580608} 
  & \displaystyle -\frac{125   e^7}{224} + \frac{7375   e^9}{2688} - \frac{120625   e^{11}}{24192} \\

7 & \displaystyle -\frac{725   e^{10}}{41472} + \frac{4375   e^{12}}{152064}
  & \displaystyle -\frac{4675   e^9}{145152} + \frac{139175   e^{11}}{580608}
  & \displaystyle \frac{1325   e^8}{8064} - \frac{111425   e^{10}}{145152} + \frac{517775   e^{12}}{387072}
  & \displaystyle -\frac{325   e^7}{4032} - \frac{975   e^9}{1792} + \frac{1148375   e^{11}}{290304} \\

8 & \displaystyle \frac{935   e^{10}}{72576} - \frac{33625   e^{12}}{399168}
  &- \displaystyle \frac{8815   e^9}{145152} + \frac{147715   e^{11}}{580608}
  & \displaystyle \frac{115   e^8}{2016} + \frac{535   e^{10}}{24192} - \frac{453545   e^{12}}{580608}
  & \displaystyle \frac{865   e^7}{4032} - \frac{29795   e^9}{16128} + \frac{801025   e^{11}}{145152} \\
\bottomrule
  \end{array}
  $
}
\caption{Hansen coefficients $X_{10}^{n,m} $ for  $0 \leq n \leq 8$, $m=0,1,2,3$ at order $12$ in $e$.}
\label{res10}
\end{table}
\FloatBarrier

\nl
Given the above tables, thanks to expansion in \ref{HansenExpansion}, one can obtain the Fourier coefficients for a fixed order of truncation in $(a,e)$. We provide some numerical example in Appendix \ref{appendicite}.

\section{Zeros of the Fourier coefficients}
Thanks to the expansion of the perturbing function $F(a,e,\ell,g)$ (see  \ref{HansenExpansion}) and to the computational power developed for Hansen coefficients shown in section \ref{HANSEN}, we can study the presence of zeros for the Fourier coefficients $f_{m,k} (a,e)$. Since the regime of application of KAM theory is for small values of $(a,e)$, we will analyze those zeros only in the region $R=(0,1)\times(0,1)$.

\nl
Given the expansion
$$
f_{m,k} (a,e)=t_{m,k} e^{|k-m|}a^{m^*} [1+g_{m,k}(a,e)]; \ \  \mbox{with} \ \ g_{m,k}(a,e)=\mathcal{O}(e^2;a)
$$
for $(m,k)\neq (0,0)$, then looking for zeros of $f_{m,k}$ is equivalent to study the presence of $(a,e) \in R$ such that
$$
g_{m,k}(a,e)=-1.
$$

\nl
For each Fourier modes $(m,k) \in \mathcal{G}^2$ we define
$$
\Gamma_{m,k} :=\{(a,e)\in R : f_{m,k}(a,e)=0 \}.
$$
We will be interested in the set of points in which both $f_{m,k}$ and $f_{2m,2k}$ or eventually $f_{m,k}$, $f_{2m,2k}$ and $f_{3m,3k}$ have a common zero. For this reason we define
$$
D_{m,k} := \Gamma_{m,k} \cap \Gamma_{2m,2k} \quad \mbox{and} \quad T_{m,k} := \Gamma_{m,k} \cap \Gamma_{2m,2k} \cap \Gamma_{3m,3k}.
$$
Moreover, let
$$
\Gamma = \bigcup_{(m,k) \in \mathcal{G}^2} \Gamma_{m,k}; \quad D = \bigcup_{(m,k) \in \mathcal{G}^2} D_{m,k}; \quad T = \bigcup_{(m,k) \in \mathcal{G}^2} T_{m,k}.
$$

\nl
Our main goal would be to find a small region in $R$ such that at least one of the above sets ($\Gamma$,$D$ or $T$) is empty. For this reason we will numerically check if those sets are empty uniformly at every order of truncation. 
\medskip

Firstly we show the graphs of the set of curves $\Gamma$. In order to do this, we have to truncate the coefficient to a certain power of $a,e$.  Thanks to the numerical power developed, we can show the curves of zeros in $R$ at 4 different orders of truncation: 5,10,20,60.

\begin{figure}[htbp] 
\centering
\begin{tabular}{cc}
    \includegraphics[width=0.45\textwidth]{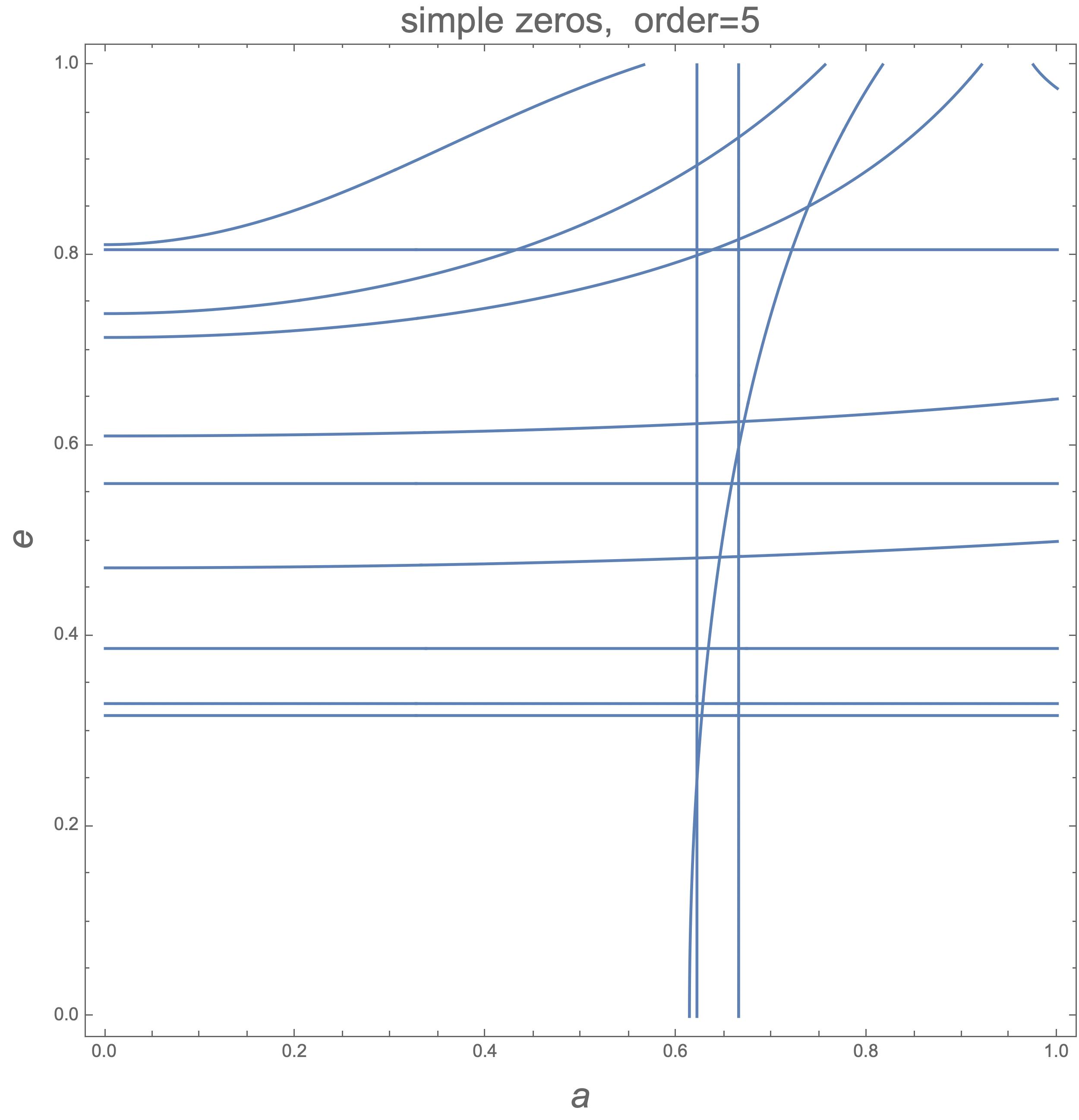} &
\includegraphics[width=0.45\textwidth]{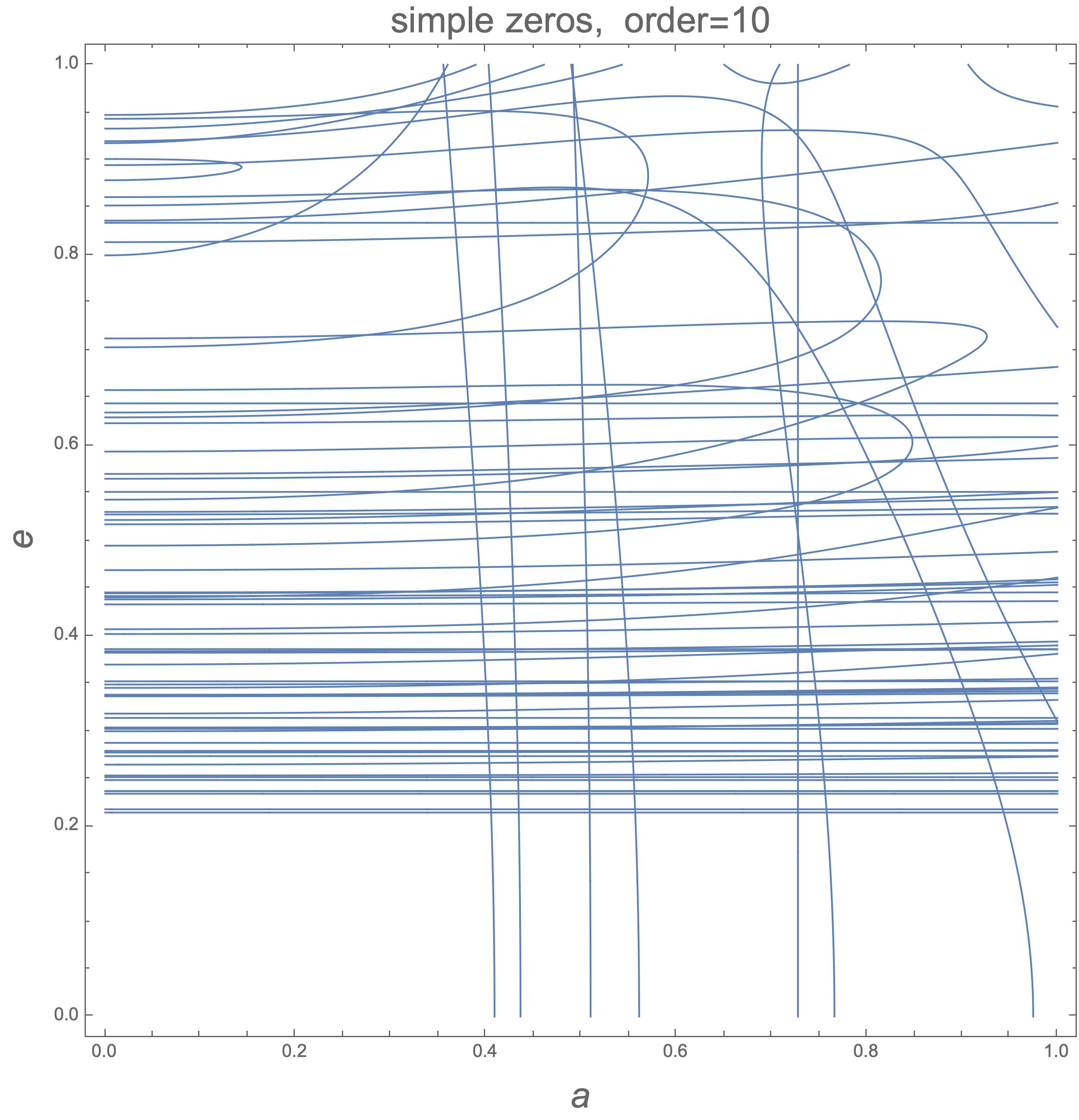} \\
    \includegraphics[width=0.45\textwidth]{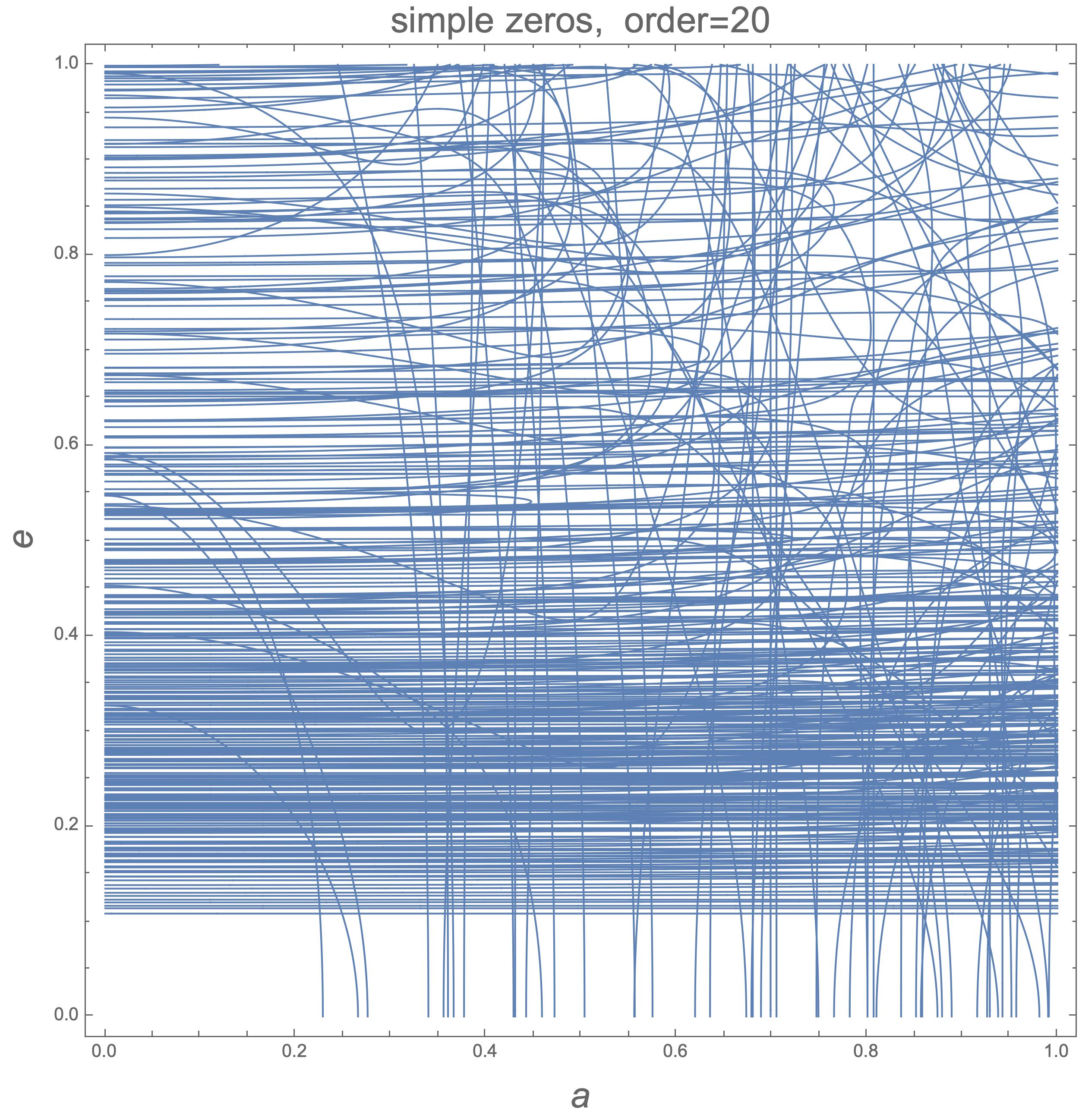} &
    \includegraphics[width=0.45\textwidth]{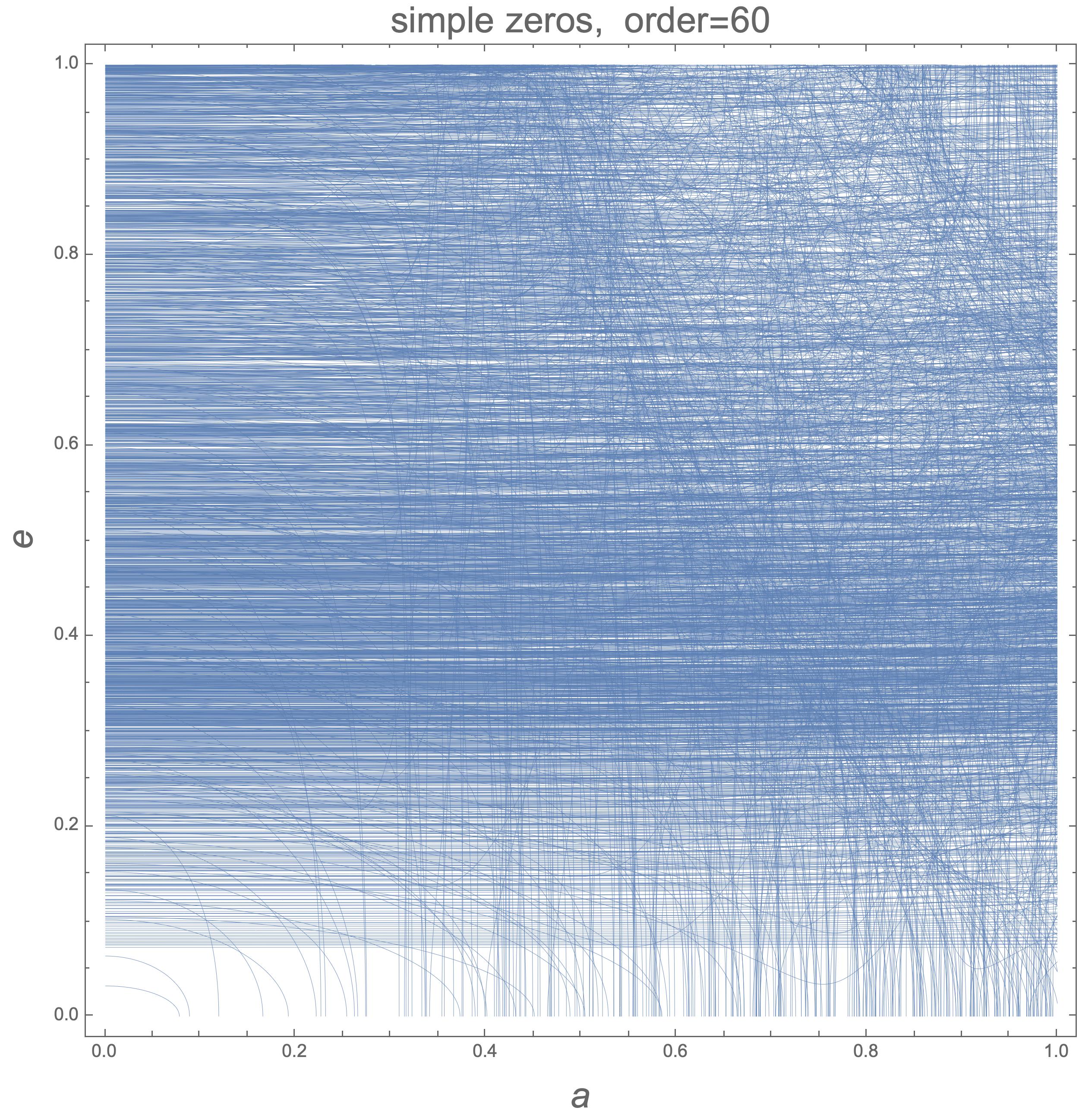} \\
\end{tabular}
\caption{The set $\Gamma$ varying the order of truncation: 5, 10, 20, 60.}
\label{ZERISINGOLI}
\end{figure}
\FloatBarrier

The results are clear: increasing the order of truncation the number of zero curves increases and their distance from the origin decline towards zero. It is impossible to fix a rectangle in which $|f_{m,k}|>0$ uniformly at every order.
\medskip

Passing to the set $D$ the situation becomes more complicated. In fact, when we consider the double coefficient, the Fourier modes start becoming bigger and since the leading powers of $a,e$ depends crucially on the values of the modes, a higher truncation is needed.

\begin{figure}[h!]
\centering
\includegraphics [scale=.5]{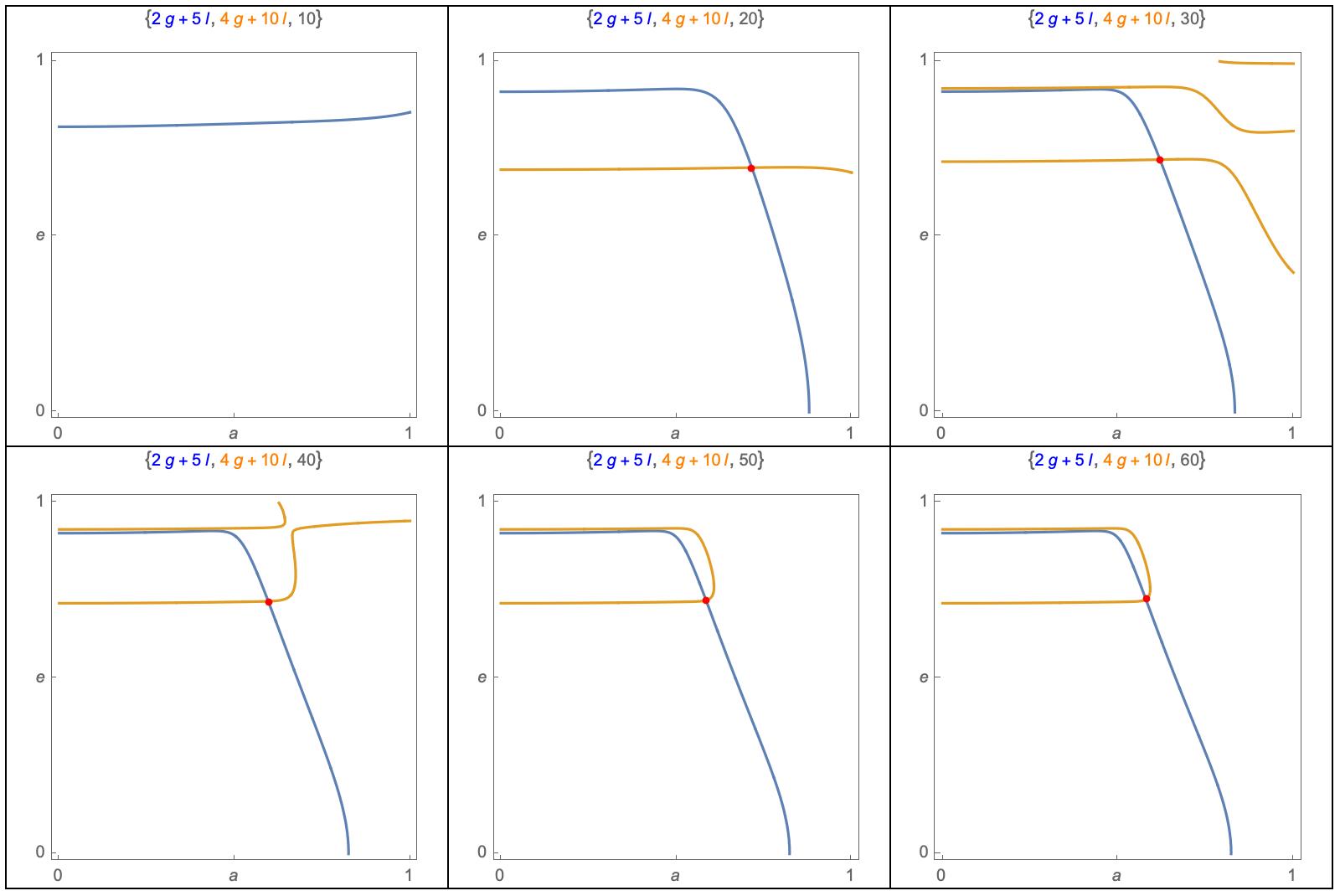}
\caption{Double zeros in the case of $(m,k)=(2,5)$ at order 10, 20, ..., 60. The red dot is the intersection point between the curves $f_{2,5}(a,e)=0$ (blue) and $f_{4,10}(a,e)=0$ (orange). Note that the curves seems to stabilize while the order tends to $60$.}
\label{Truncationdouble}
\end{figure}

To show the above fact, we present in figure \ref{Truncationdouble} how the order of truncation (from 10 to 60) changes the graphs of the double case for $(m,k)=(2,5)$. 
Note that the curves change a lot at low orders but seems to stabilize while the order tends to $60$

\FloatBarrier

\begin{figure}[htbp] 
\centering
\begin{tabular}{cc}
    \includegraphics[width=0.45\textwidth]{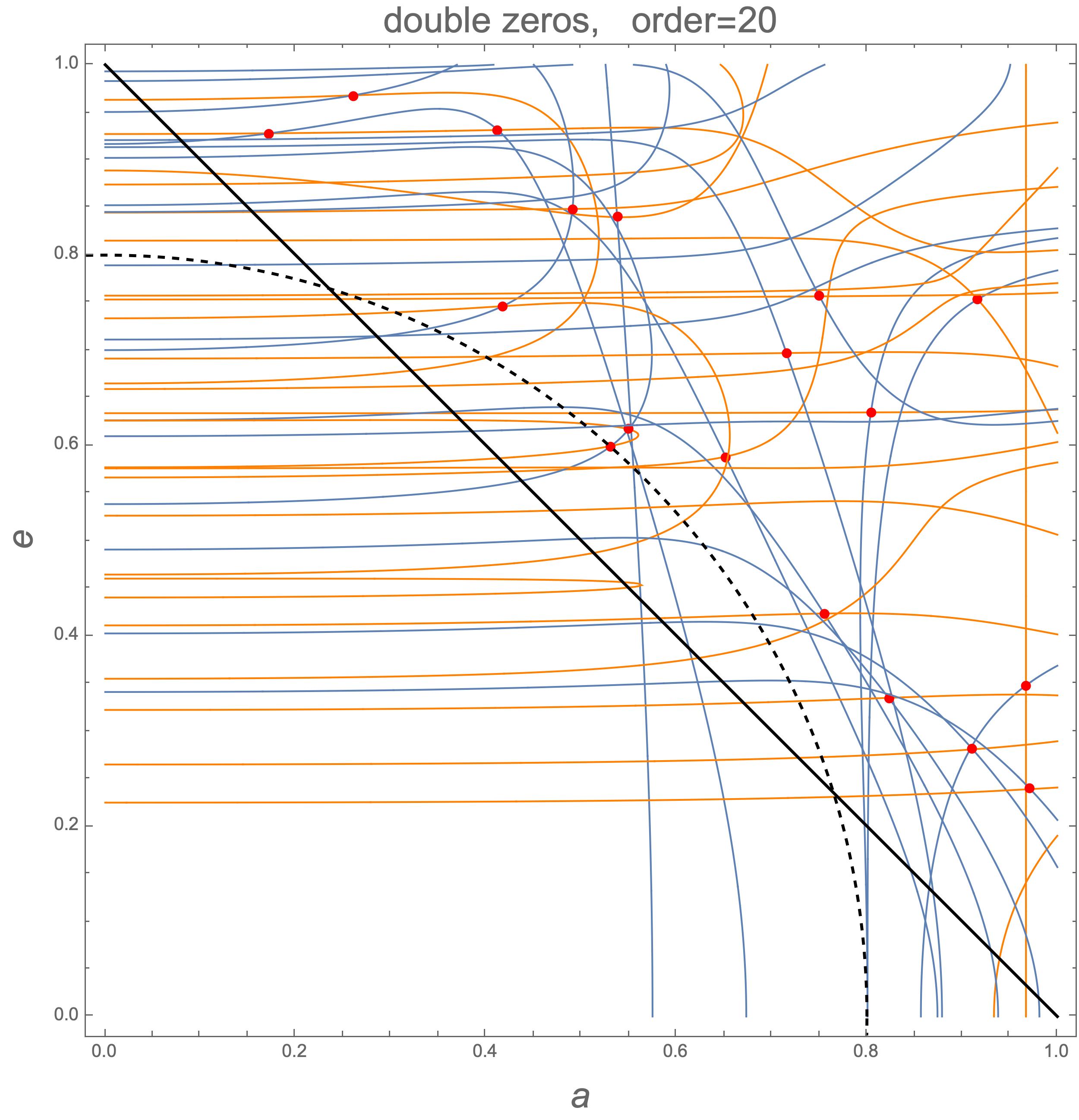}  &
\includegraphics[width=0.45\textwidth]{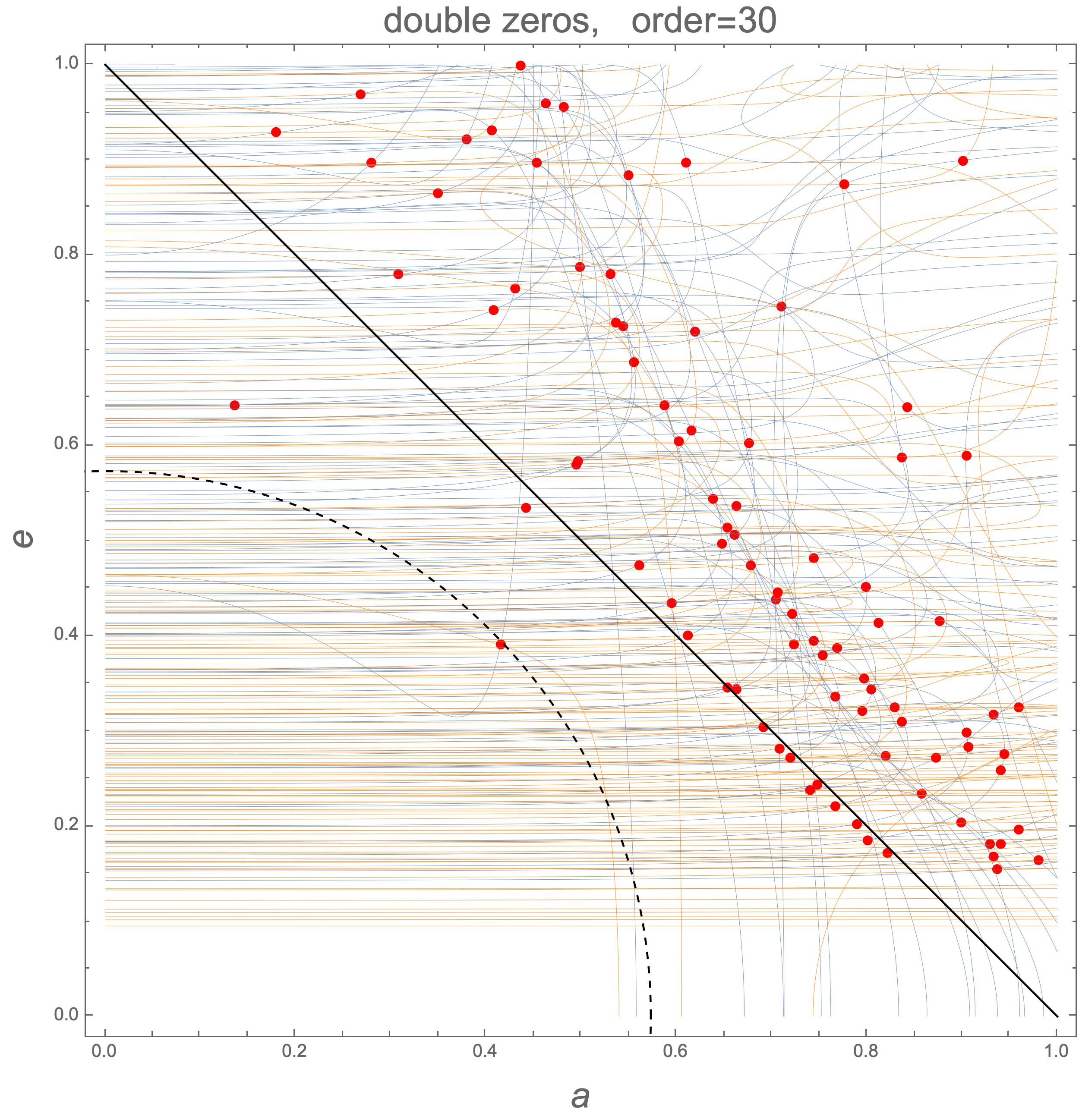}  \\
    \includegraphics[width=0.45\textwidth]{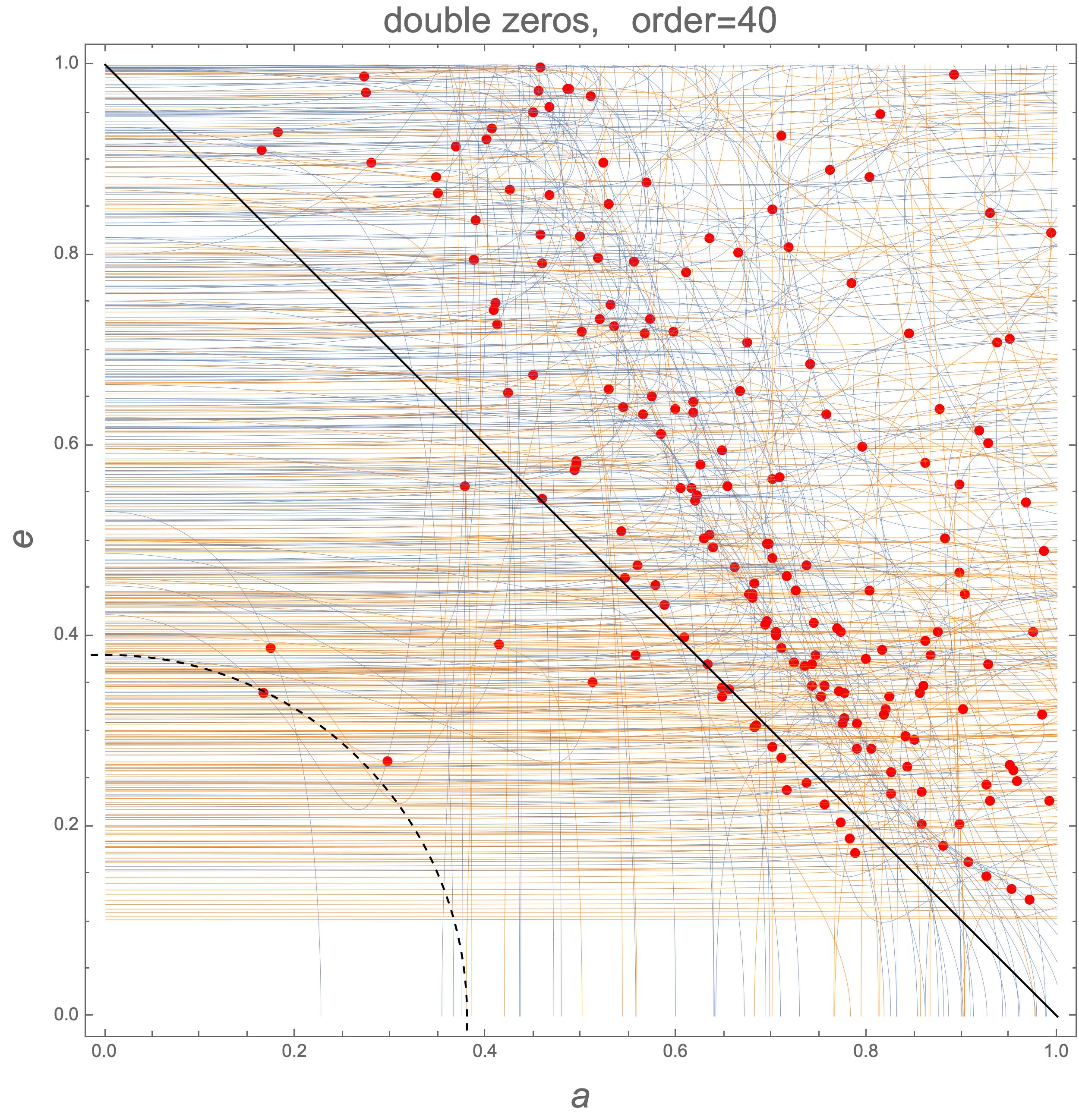}   &
    \includegraphics[width=0.45\textwidth]{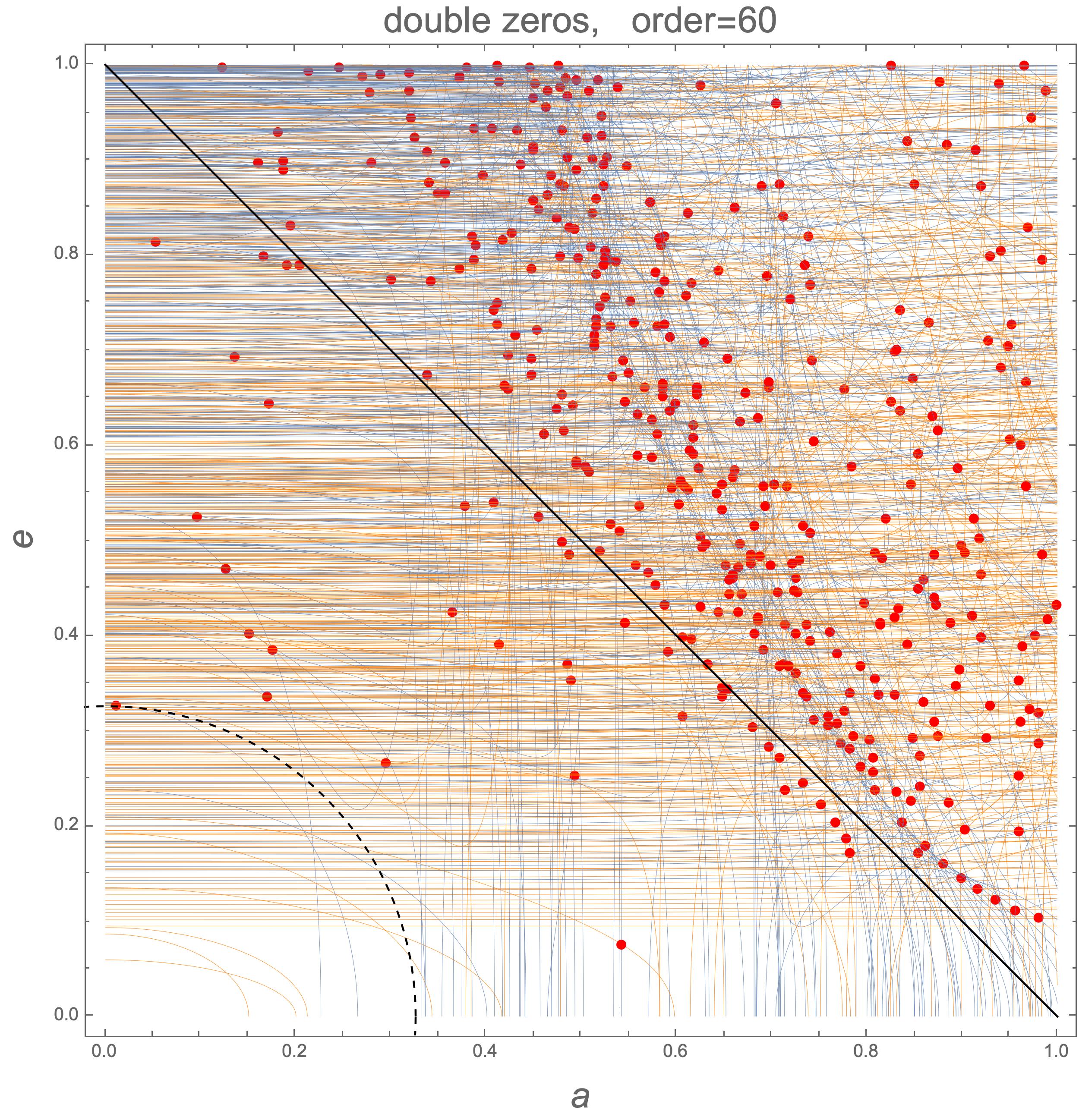} \\
\end{tabular}
\caption{The set $D$ varying the order of truncation (20,30,40,60). Each red dot is an intersection point between two curves: $f_{m,k}(a,e)=0$ (blue) and $f_{2m,2k}(a,e)=0$ (orange) for $(m,k) \in \mathcal{G}^2$. The radius of the dashed circle indicates the distance of the set $D$ from the origin.}
\label{ZERIDOPPI}
\end{figure}

\begin{figure}[h!]
\centering
\includegraphics [scale=.42]{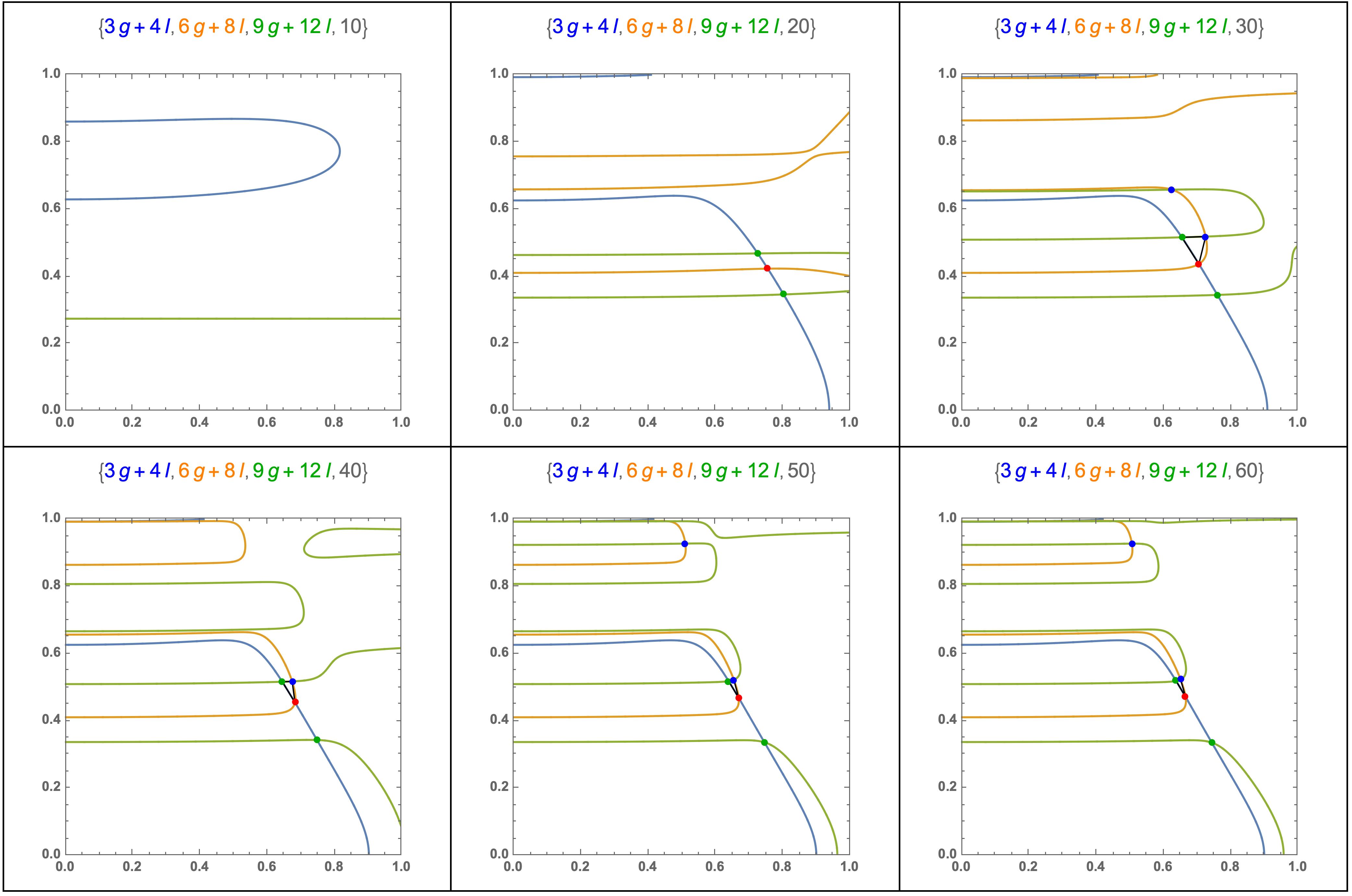}
\caption{Zero curves $f_{3,4}(a,e)=0$ (blue), $f_{6,8}(a,e)=0$ (orange) and $f_{9,12}(a,e)=0$ (green), with their interection points at order 10, 20, ..., 60. Note the triangle formed by the three curves starting at order 30.}
\label{truncationtriple}
\end{figure}

We then compute all the points belonging to the set $D$ up to order 60 and show these points at four different order of truncation (20, 30, 40, 60), together with the zero curves that generate them (see figure \ref{ZERIDOPPI}).

\FloatBarrier

Also in figure \ref{ZERIDOPPI} we can see how, increasing the order, the number of points in $D$ increases and they come closer and closer to the origin. So we conclude that it is not possible to have a region $R'\subset R$ such that $R' \cap D=\emptyset$ uniformly in the order of truncation.
\bigskip

\nl
Finally we try to do the same analysis for the set $T$. Note that for each $(m,k)$ we now check also higher modes $(3m,3k)$, and the corresponding zero curves do not always stabilize up to order 60. 

Here we report an example for the case $(m,k)=(3,4) \in \mathcal{G}^2$. At sufficiently high order, double zeros come closer and a curved triangle is formed by the three curves $f_{3,4}(a,e)=0$, $f_{6,8}(a,e)=0$, $f_{9,12}(a,e)=0$. It appears to shrink at higher orders.

\FloatBarrier

\begin{figure}[h!]
\centering
\includegraphics [scale=.4]{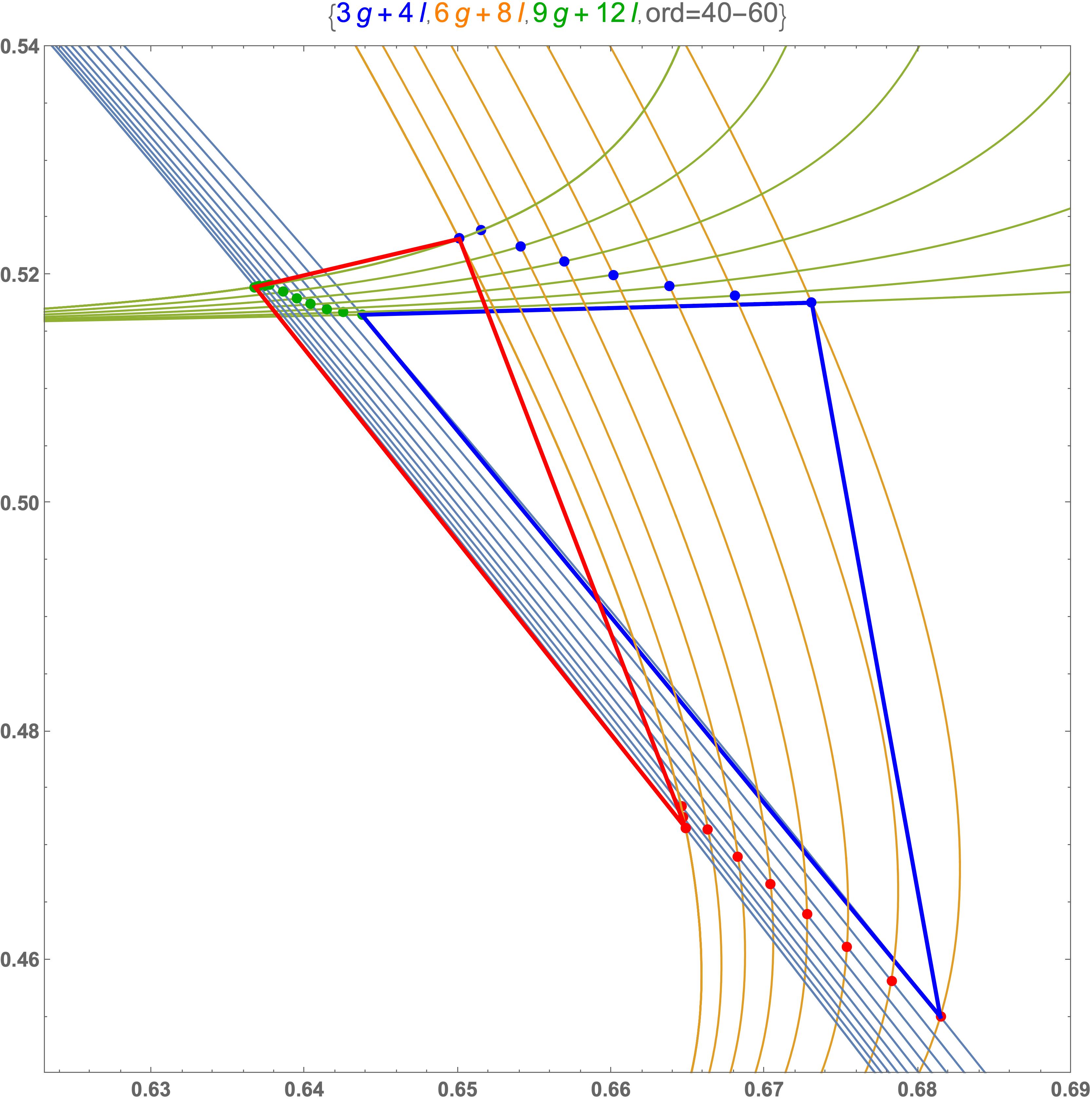}
\caption{Zoom of figure \ref{truncationtriple} for orders 40, 42, ... 60. Zero dots are red for modes $(3,4)$ and $(6,8)$, green for $(3,4)$ and $(9,12)$, blue for $(6,8)$ and $(9,12)$. The area of the blue triangle (order 40) is $9.22\cdot 10^{-4}$ while the red area (order 60) is $3.75\cdot 10^{-4}$. Note the convergences as the order increases.}
\label{zerizoom}
\end{figure}

For this close to triple intersection, we can now zoom in at a very small scale and see that, varying the order, the triangle seems to stabilize (see figure \ref{zerizoom}).

\FloatBarrier

Now, to compute the set $T$ we first have to check all possible intersection points between $f_{m,k}$, $f_{2m,2k}$ and $f_{3m,3k}$ varying $(m,k) \in \mathcal{G}^2$ at different order of truncation. Up to order 60 there are no points of ‘‘triple’’ intersection (i.e. $T = \emptyset$). Since it would be impossible to show all possible cases, we report here only the three curves of zeros ($\{ (a,e) \in R : f_{jm,jk}(a,e)=0\}$ for $j=1,2,3$) which form triangles of minimal area less than $10^{-4}$.


\medskip

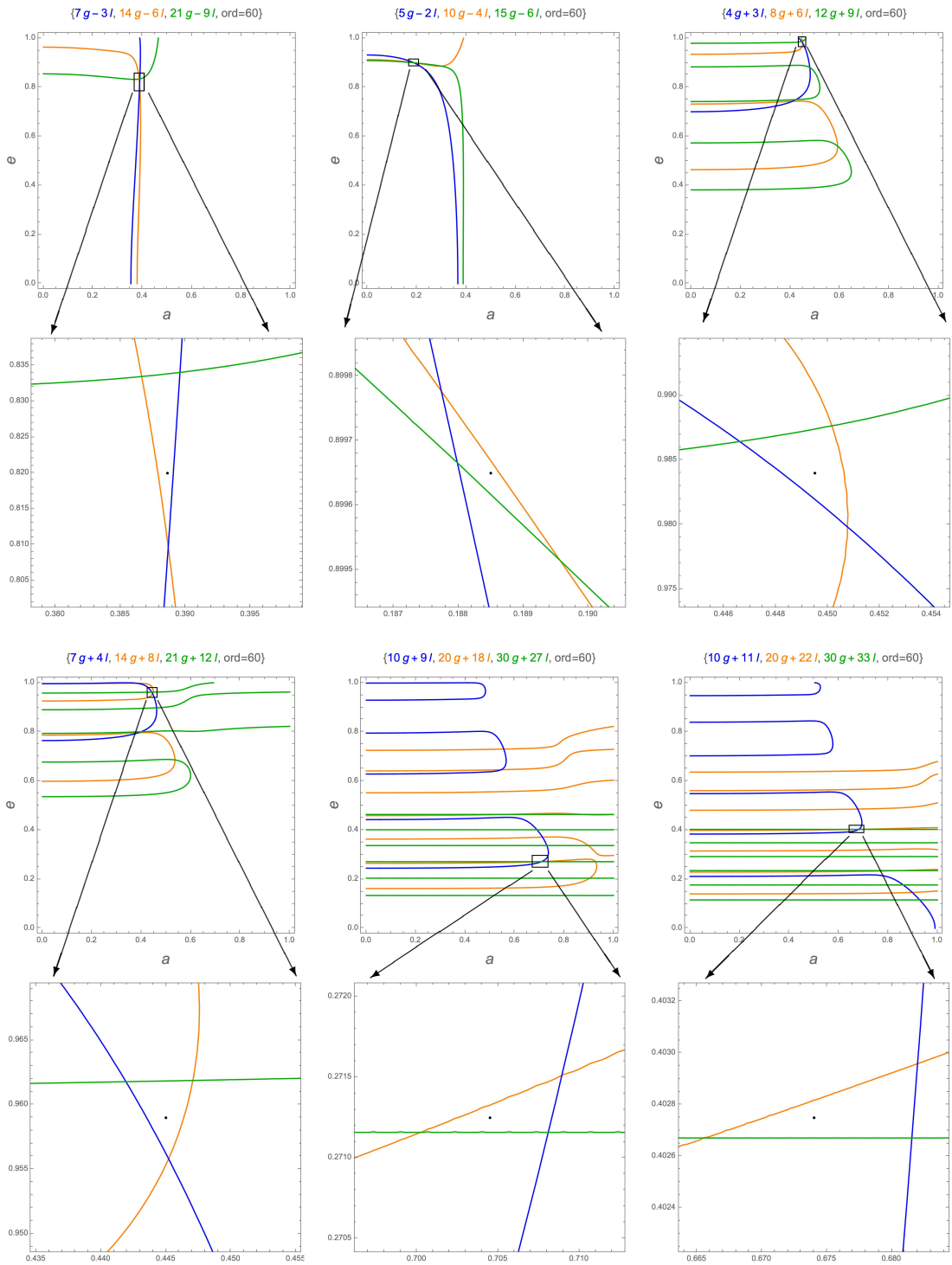
\begin{figure}[htbp] \label{figuresix}
\centering
\begin{tabular}{c}
    \includegraphics[width=0.95 \textwidth]{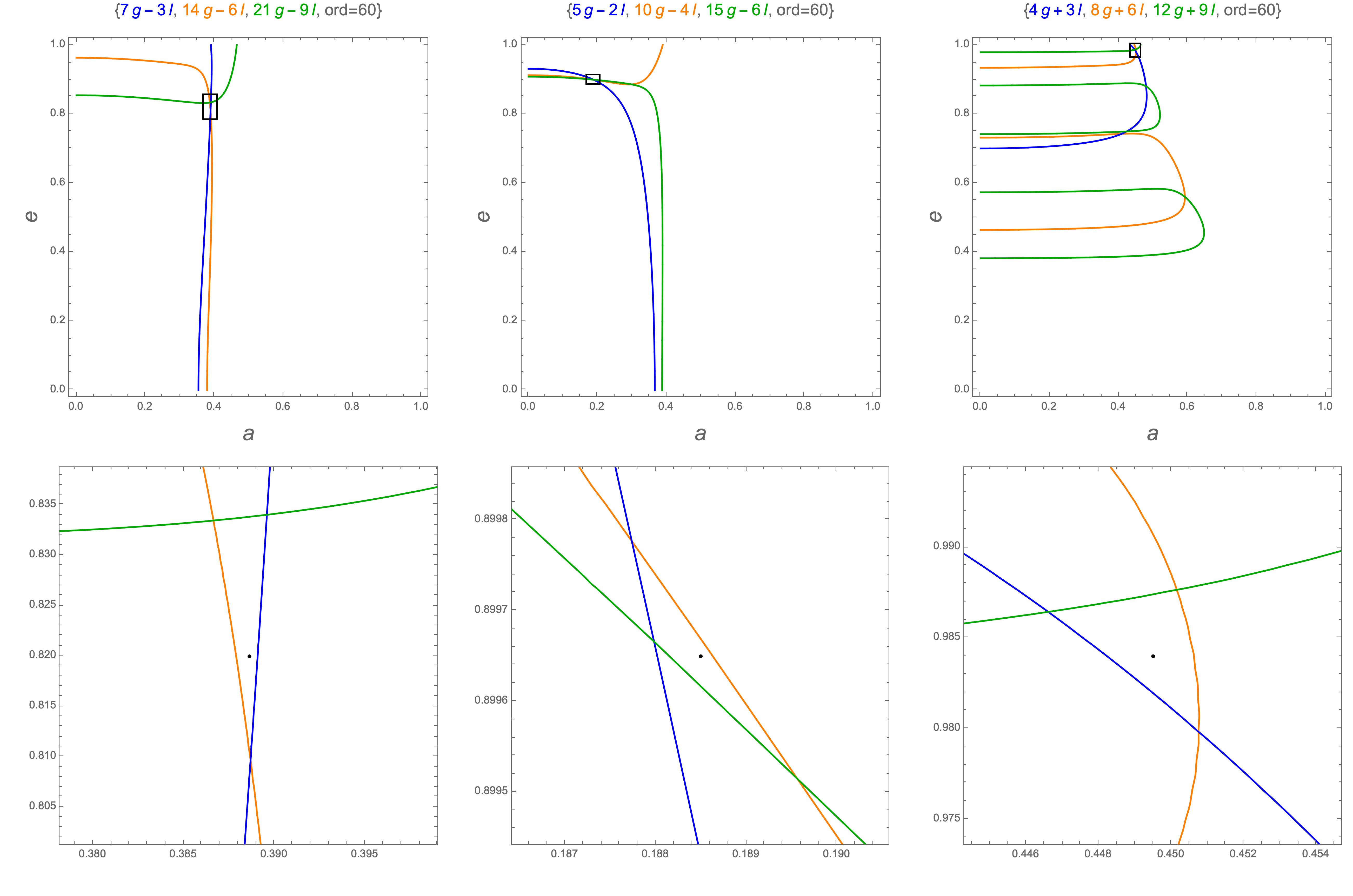}
    \put(-350,240){\vector(1,-2){52}}  \put(-357,240){\vector(-1,-3){35}}
    \put(-230,250){\vector(2,-3){76}}  \put(-237,250){\vector(-1,-4){28}}
    \put(-65,260){\vector(1,-2){60}}  \put(-70,260){\vector(-1,-3){40}}\\
    \includegraphics[width=0.95 \textwidth]{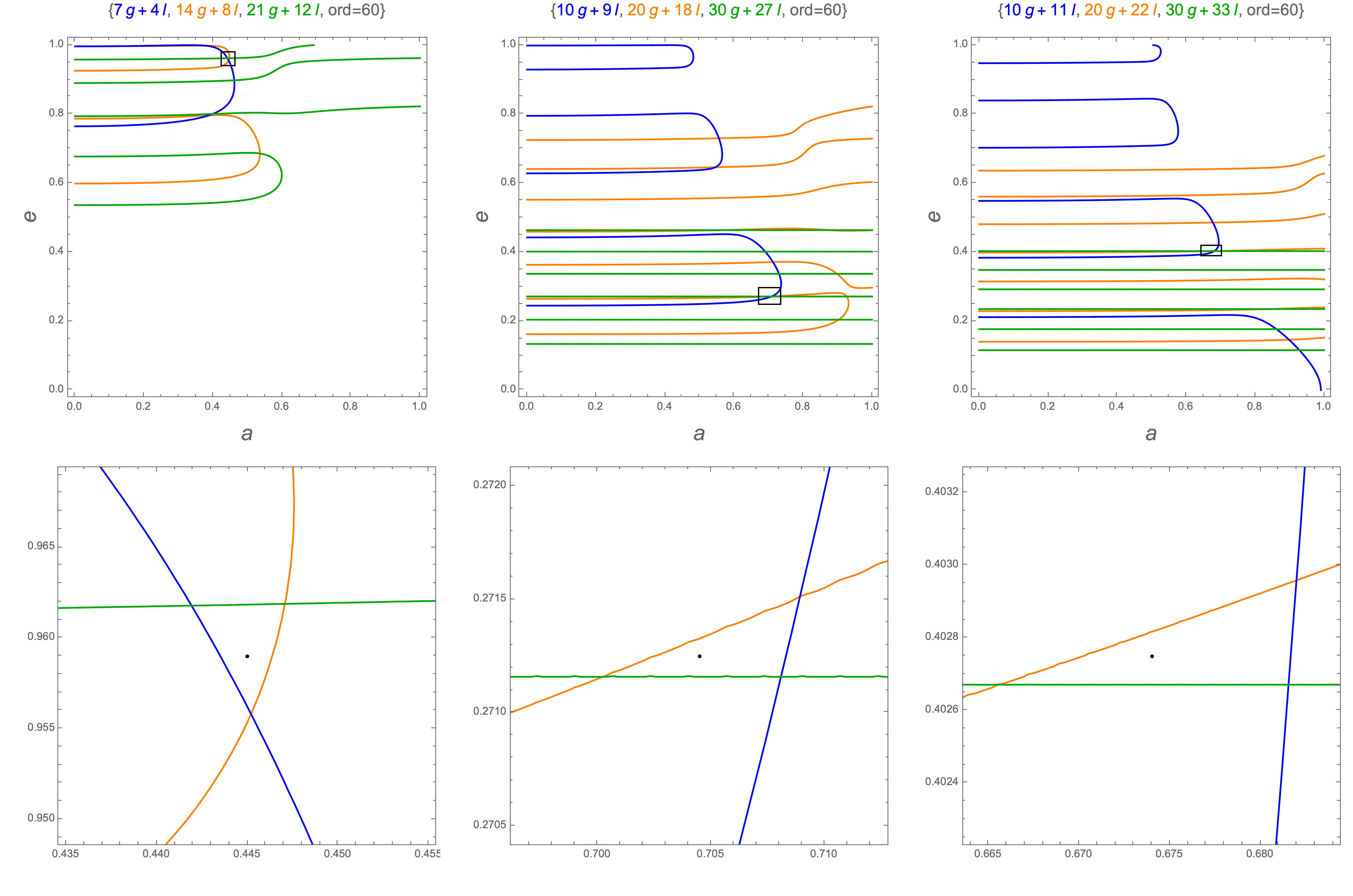}
    \put(-346,256){\vector(1,-2){60}}  \put(-351,256){\vector(-1,-3){40}}
    \put(-177,182){\vector(2,-3){31}}  \put(-184,182){\vector(-3,-2){70}}
    \put(-41,197){\vector(1,-2){31}}  \put(-47,197){\vector(-1,-1){62}}\\
\end{tabular}
\caption{The only six cases (first and third row) in which the three zero curves ($\{ (a,e) \in R : f_{jm,jk}(a,e)=0\}$, $j=1,2,3$) form triangles (second and fourth row) of minimal area less than $10^{-4}$.}
\label{zeri3}
\end{figure}
\FloatBarrier

From the above figure, we emphasize that the most interesting case is the mode $(m,k)=(5,-2)$ for which the zero curves form a triangle of area less than $10^{-7}$, incenter $(a,e)=(0.18799,0.89970)$ and radius of the inscribed circle $3.78\ 10^{-5}$. We zoom in (see Figure 7) and check how, varying the order from 58 to 60, this area changes. The figure does not stabilize at truncation order 60, since the expansion depends strongly on how big is the value of $3m$ and $3|m-k|$ as we can see from section \ref{section2}.

\begin{figure}[h!]
\centering
\includegraphics [scale=.5]{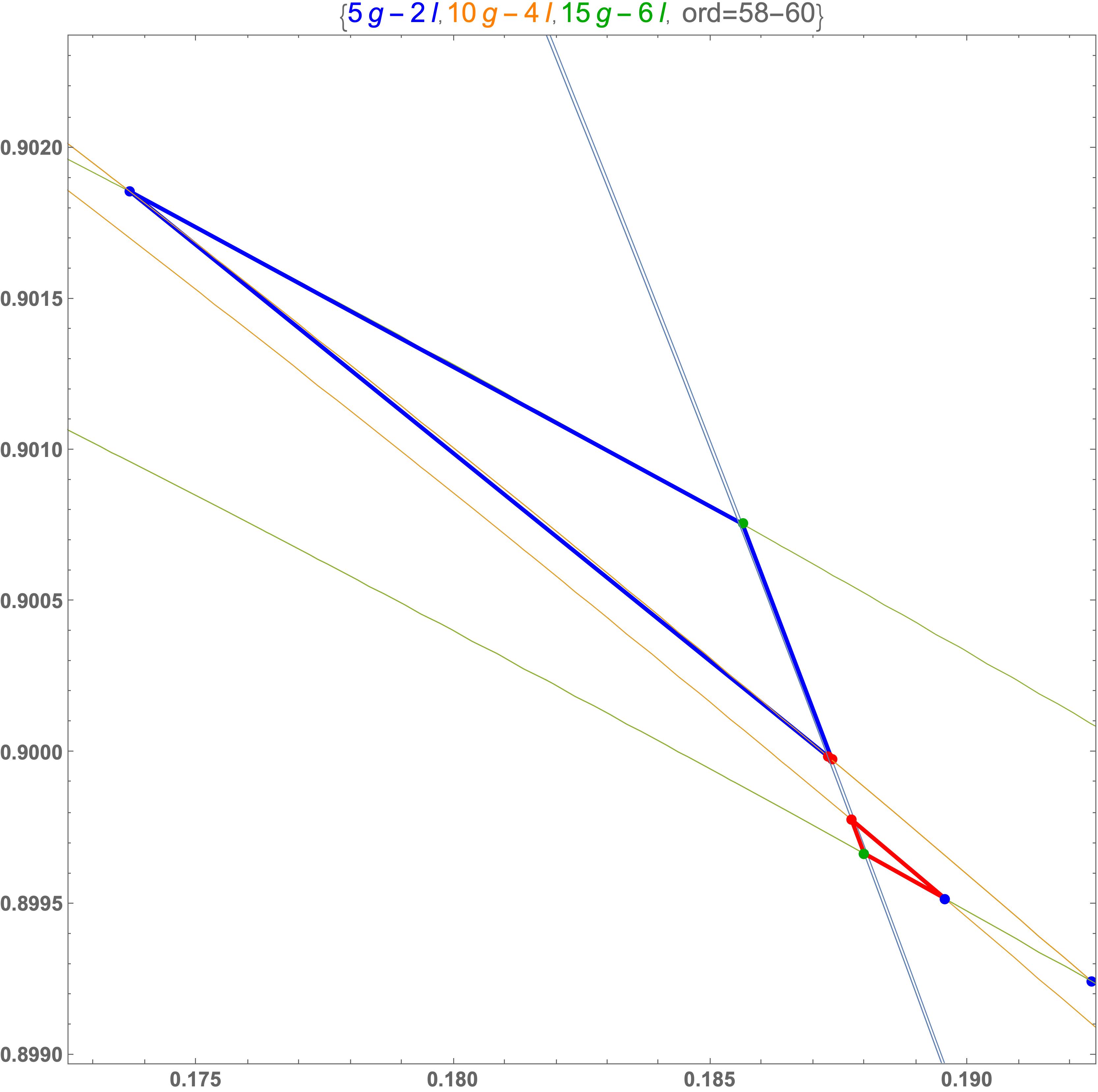}
\caption{Zoom on the close to triple intersection case for $(m,k)=(5,-2)$ (see Figure 6). The blue triangle, at order $58$, has area  $3.68\ 10^{-6}$ and radius of the inscribed circle $2.67\ 10^{-4}$. The red triangle, at order $60$, has area  $6.97\ 10^{-8}$ and radius of the inscribed circle $3.78\ 10^{-5}$. }
\label{zerizoom52}
\end{figure}
\FloatBarrier

So the triangle does not stabilize at order 60 (see figure \ref{zerizoom52}) and we cannot exclude the possibility that a triple zero of $f_{5,-2}(a,e), f_{10,-4}(a,e)$ and $f_{15,-6}(a,e)$ could appear by increasing the order of truncation.

\nl
Summing up, for the case of PCR3BP perturbing function, there is no hope to find a small region with empy sets $\Gamma$ and $D$ (figures \ref{ZERISINGOLI} and \ref{ZERIDOPPI}) , but probably the case with 3 or more zeros of coefficients could work (figure \ref{zeri3}). 
This suggests that in order to hope for physical applications of the estimates on secondary tori, we should certainly develop and understand in a deeper way, even from a theoretical point of view,  what is suggested in Barbieri-Niedermann \cite{TesiSanti}, developing KAM for secondary tori even in the case of action--dependent perturbation using quantitative Morse-Sard theory.

\appendix


\small

\section{Tables of truncated Fourier coefficients }\label{appendicite}

We report here some tables showing the Fourier coefficients of the perturbing function truncated at order 15 in power of $a,e$.
This is due to the expansion in \ref{HansenExpansion} togheter with the numerical power developed in section 3 for Hansen coefficients.

\begin{figure}[h!]
\centering
\includegraphics [scale=.53]{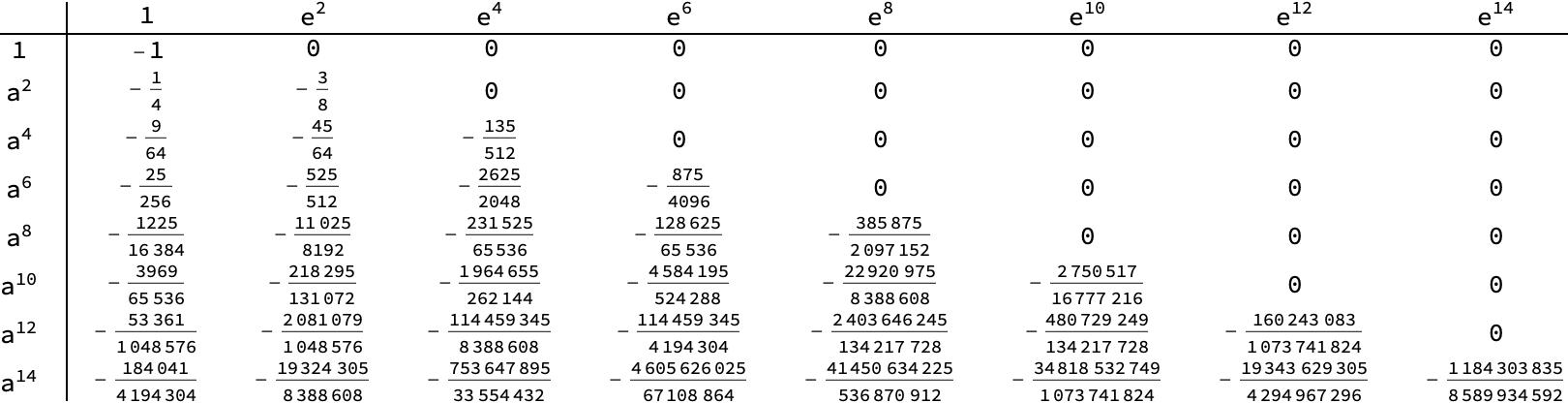}
\caption{Coefficients of $f_{0,0}(a,e)$ at order 15 in $a$ and $e$.}
\label{tab00o15}
\end{figure}
\FloatBarrier

\begin{figure}[h!]
\centering
\includegraphics [scale=.51]{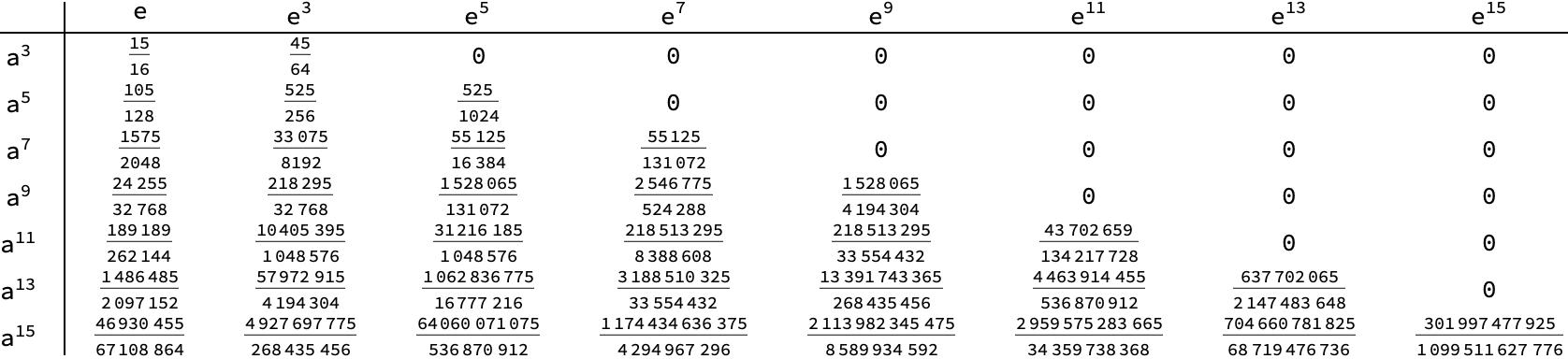}
\caption{Coefficients of $f_{1,0}(a,e)$ at order 15 in $a$ and $e$.}
\label{tab10o15}
\end{figure}
\FloatBarrier

\begin{figure}[h!]
\centering
\includegraphics [scale=.44]{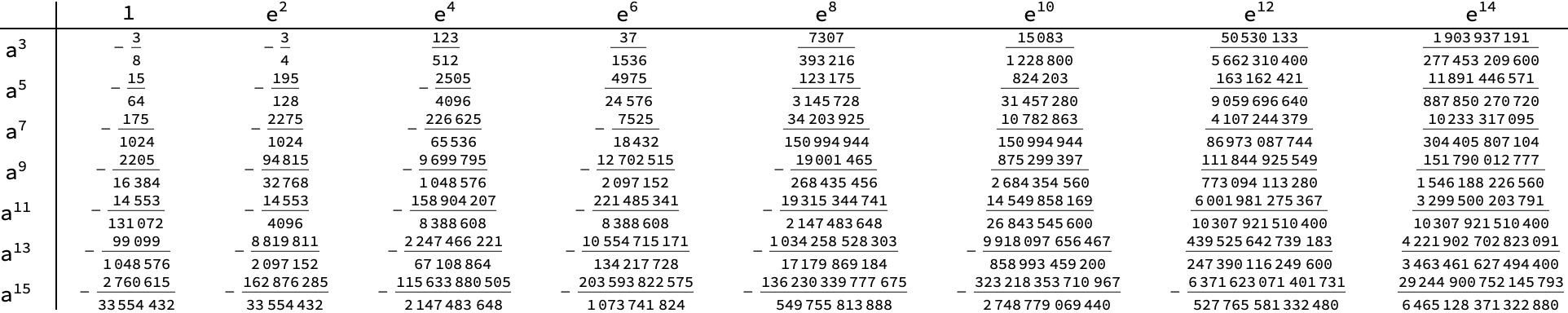}
\caption{Non-zero coefficients of $f_{1,1}(a,e)$ at order 15 in $a$ and $e$.}
\label{tab11o15}
\end{figure}
\FloatBarrier

\begin{figure}[h!]
\centering
\includegraphics [scale=.54]{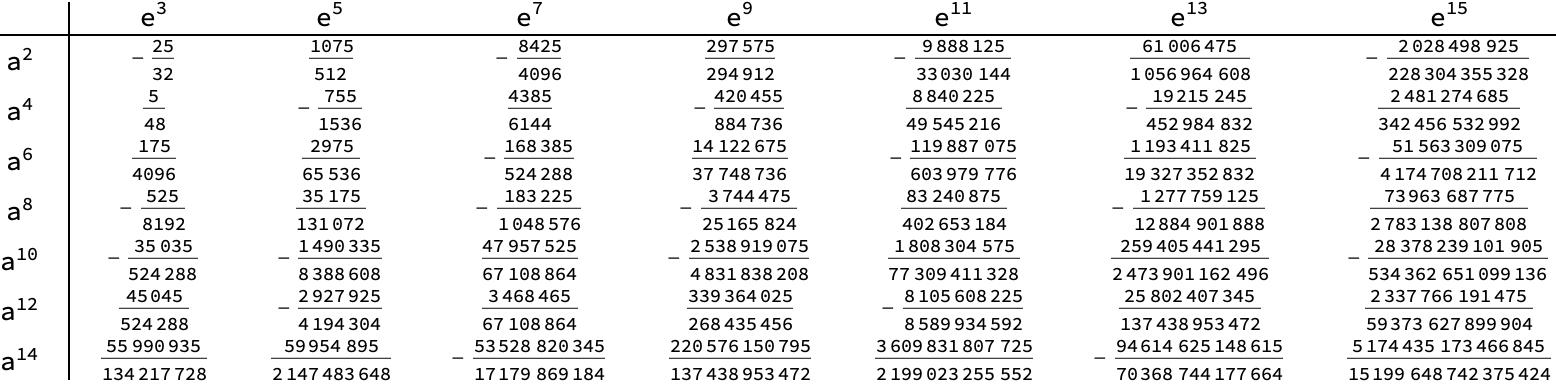}
\caption{Non-zero coefficients of $f_{2,5}(a,e)$ at order 15 in $a$ and $e$.}
\label{tab25o15}
\end{figure}
\FloatBarrier

\begin{figure}[h!]
\centering
\includegraphics [scale=.5]{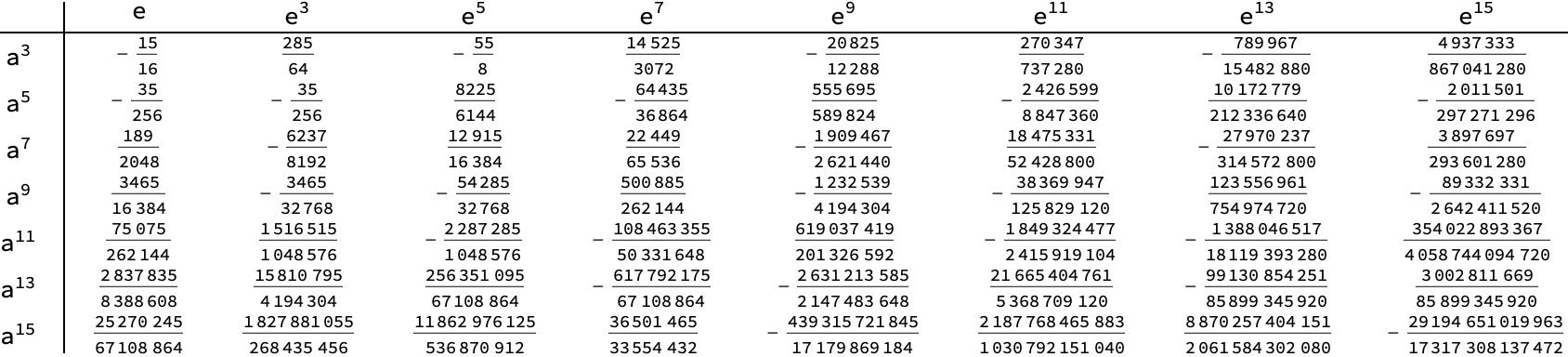}
\caption{Non-zero coefficients of $f_{3,4}(a,e)$ at order 15 in $a$ and $e$.}
\label{tab34o15}
\end{figure}
\FloatBarrier

\newpage

\facciatabianca

\end{document}